\documentclass[11pt]{amsart}
\usepackage{pstricks-add}
\usepackage{amsfonts}
\usepackage{amssymb, arcs}
\usepackage{amsxtra}
\usepackage{amsthm}
\usepackage{amsfonts, amscd, pstricks}
\usepackage{graphicx}
\usepackage{asymptote}
\usepackage{physics}
\include{std}
\usepackage[all]{xy}

\newdimen\minCDarrowwidth
\usepackage{amsmath}
\usepackage{amsfonts}
\usepackage{amssymb}
\usepackage{amscd}
\usepackage{verbatim}

\newcommand{\marginlabel}[1]%
 {\mbox{}\marginpar{\raggedleft\hspace{0pt}\bfseries\sf#1}}

\def\AA{{\mathbb A}}
\def\CC{{\mathbb C}}
\def\PP{{\mathbb P}}
\def\ZZ{{\mathbb Z}}

\def\cT{\mathcal{T}}
\def\cO{\mathcal{O}}
\def\cM{\mathcal{M}}

\def\cH{\mathcal{H}}

\def\ra{\rightarrow}
\def\dra{\dashrightarrow}

\def\ol{\overline}

\newtheorem{lemma}{Lemma}[section]
\newtheorem{theorem}[lemma]{Theorem}
\newtheorem{corollary}[lemma]{Corollary}
\newtheorem{proposition}[lemma]{Proposition}

\theoremstyle{definition}

\newtheorem{example}[lemma]{Example}

\newtheorem{notation}{Notation}

\numberwithin{equation}{section}

\newcommand{\bean}{\begin{eqnarray}}
\newcommand{\eean}{\end{eqnarray}}
\newcommand{\bl}{\begin{lemma}}
\newcommand{\el}{\end{lemma}}
\newcommand{\br}{\begin{red}}
\newcommand{\er}{\end{red}}
\newcommand{\bp}{\begin{proposition}}
\newcommand{\ep}{\end{proposition}}
\newcommand{\bt}{\begin{theorem}}
\newcommand{\et}{\end{theorem}}
\newcommand{\bpr}{\begin{proof}}
\newcommand{\epr}{\end{proof}}
\newcommand{\bea}{\begin{eqnarray*}}
\newcommand{\eea}{\end{eqnarray*}}
\newcommand{\be}{\begin{displaymath}}
\newcommand{\ee}{\end{displaymath}}

\begin{document}

\title{Singular plane sections and the conics in the Fermat quintic threefold}

\author[Anca Musta\c{t}\v{a}]{Anca~Musta\c{t}\v{a}}
\address{School of Mathematical Sciences, University College Cork, Ireland }
\email{{\tt a.mustata@ucc.ie}}

\date{\today}

%\subjclass{14J30 (Primary) 14J32 (Secondary)}
%\keywords{Hilbert scheme, quintic threefold, conics}

\begin{abstract}

We present explicit equations for the space of conics in the Fermat quintic threefold $X$, working within the space of plane sections of $X$ with two singular marked points. This space of two-pointed singular plane sections has a  birational morphism to the space of bitangent lines to the Fermat quintic threefold, which in its turn is birational to a 625-to-1 cover of $\PP^4.$ We illustrate the use of the resulting equations in identifying special cases of one-dimensional families of conics in $X.$
\end{abstract}

\maketitle
\begin{center}
 \emph{ Dedicated to Herb Clemens, to whom I am grateful for his mathematical insights, generosity, patience, encouragement and empathy. }   
\end{center}

\bigskip

\section{Introduction}

The Fermat quintic threefold $X$ is defined in $\PP^4$ by the polynomial  \bean  \label{fermat eq} F(x_0, x_1, x_2, x_3, x_4) = \sum_{i=0}^{i=4}x_i^5. \eean
As such it has a large group of automorphisms isomorphic to $S_5\times \ZZ_5^4$, which is at the root of many wonderful properties of this manifold. Some of these properties are shared with, and gain a special significance when looking at the Fermat quintic as a member of the Dwork pencil $X_\psi \hookrightarrow \PP^4$, given by defining polynomials
 \[F_{\psi} = \sum_{i=0}^{i=4}x_i^5-5\psi\prod_{i=0}^{i=4}x_i. \]
 
Such is the case of the spaces of rational curves in degrees 1 and 2 lying on these threefolds. 
It is well known that a generic quintic threefold contains 2875 lines and 609,250 conics. On the other hand, all members of the Dwork pencil contain continuous families of lines and of smooth conics. 

The lines in the Dwork pencil have been studied in \cite{ak1}, \cite{ak2}, \cite{eu}, ending with the remarkably precise and beautiful presentation in \cite{candelas lines} and its interpretation in  \cite{zagier}. We briefly outline the history here. The hyperplane $V(x_0+x_1)$ intersects the Fermat quintic $X$ at a cone of vertex $[1:-1:0:0:0]$ over the Fermat quintic curve 
\bea x_2^5+x_3^5+x_4^5=0. \eea
A total of 50 such cones can be obtained by  permutations of coordinates along with the action of the group $\ZZ_5^4$ by
\bea (h, [x_0:x_1:x_2:x_3:x_4]) \longrightarrow  [x_0:x_1\mu^{h_1}:x_2\mu^{h_2}:x_3\mu^{h_3}:x_4\mu^{h_4}] \quad  \mbox{ for }  \quad  h=(h_i)_i\in \ZZ_5^4, \eea 
where $\mu$ is a fifth root of 1. The lines through the vertex lying in these cones contribute to 50 components of the Hilbert scheme $\cH_1(X)^{\mbox{red}}$, which are smooth Fermat curves of genus $g=6$. In \cite{ak1}, Alberto Albano and Sheldon Katz showed that these are precisely all the 1-dimensional components of $\cH_1(X)$, and they each come with multiplicity 2. In addition, there are 375 special points  of $\cH_1(X)$ with multiplicity 5, each representing an element  in the orbit of the line
\bean \label{isolated line } [u:v] \longrightarrow [u:-u: v: -v: 0] \quad  \mbox{ for }  \quad [u:v]\in\PP^1 \eean 
under the action of $S_5\times \ZZ_5^4$. Thus according with the formula by Herb Clemens and Holger Kley in \cite{ck}, the family $\cH_1(X)$ contributes to the expected
\bean \label{contribution to lines } 50 \cdot 2 \cdot (2g-2)+ 375 \cdot 5 = 2875 \eean
number of lines on a generic deformation of $X$.

Within the Dwork pencil, symmetry under an $S_3$ subgroup of $S_5$ led Bert van Geemen to 5,000 lines in each $X_\psi$. Since this number is larger than the 2875 lines on a generic quintic threefold, this implied the existence of 1-dimensional components of $\cH_1(X_\psi)$ (see \cite{ak2}). Indeed, in \cite{eu} we showed that in general $\cH_1(X_\psi)$ consists of two isomorphic curves $C_{\psi}^{\pm}$ of genus 626, together with the 375 isolated points shared with $\cH_1(X)$. Of the permutation group $S_5$ of $X_\psi$, the even permutations act on each curve $C_{\psi}^{\pm}$ separately, while the odd ones swap $C_\psi^+$ and $C_\psi^-$. 
%(This can be better observed after a change of parameter $\psi \to \varphi_\pm(\psi)$).
Each member of the family $X_\psi$ also has automorphisms given by the subgroup $\ZZ_5^3$ of $\ZZ_5^4$ consisting of tuples with zero sum, and an additional $\ZZ_5$ acts on the parameter $\psi$. A new parameter $\varphi$ can be set to account for quotient by this last $\ZZ_5$ and for the two components above. 

In \cite{candelas lines}
Philip Candelas, Xenia de la Ossa,  Bert van Geemen and  Duco van Straten found explicit parametrisations for the curves $C_\varphi^+$ and $C_\varphi^-$, which allowed them to identify the quotients $C_\varphi^\pm/\ZZ_5^3$ (up to desingularisations) with the special pencil of plane sextics known as the Wiman pencil. A first member of this family was discovered by Wiman in 1897 (\cite{wiman}), as having an unusual automorphism group isomorphic with $S_5$;
%%check this. was it $A_5$?
the other members of the pencil share the same property and were discovered by Edge in \cite{edge}.

In \cite{zagier}, Don Zagier observed that the remarkable coordinate presentation of $C_\varphi^\pm/\ZZ_5^3$ in  \cite{candelas lines} has a more natural invariant theory based - interpretation based on regarding the ambient space for these curves as $\overline{M}_{0,5}$. The moduli space of rational curves with 5 marked points $\overline{M}_{0,5}$ 
%can be constructed by blow-up of $\PP^2$ at 4 points, and it can also be thought of as a quotient of the Grassmanian $G(2,5)=G(\PP^1,\PP^4)$ by the action of the torus $(\CC^*)^4.$ It 
comes with the natural action of the permutation group $S_5.$ Invariants under this action were also used in \cite{zagier} to recover the degree in $\PP^4$ of the hypersurface  covered by all the lines in $(X_\psi)_\psi$ - this invariant was calculated in \cite{eu} by intersection of Schubert cycles on the Grassmanian $G(2,5)=G(\PP^1,\PP^4).$ 

So far in the study of lines, the Dwork pencil has been grounds for some remarkable coincidences and an interplay of methods from deformation to invariant theory, leading to a surprisingly beautiful and concrete presentation. It is natural to wonder if such features will be discernible for higher degree curves. On the one hand, a nice parametrisation will be harder to find due to the increase in the dimension of the ambient space of degree $d$ rational curves in $\PP^4$, which is $(5d+1).$  On the other hand, it is to be expected that the rich symmetry will continue to play a role.

Indeed, the context proposed by Don Zagier in \cite{zagier} can be extended to degree $d$ rational curves as follows: over suitably well-balanced open sets, the moduli space of stable maps $\overline{M}_{0,0}(\PP^4,d)$ has natural local covers $\overline{M}_{0,0}(\PP^4,d,t)$, which depend on a fixed coordinate $t$ on $\PP^1$ and are $(\CC^*)^4$--bundles over open sets in $\overline{M}_{0,5d}$. Concretely, a map
\bea  \gamma: \PP^1 & \longrightarrow & \PP^4 \\
\left[ t:1\right] & \longrightarrow & [c_i\prod_{j=1}^d(t-r_{ij})]_{i\in\{0,1,2,3,4\}},\eea 
has an associated tuple $(\mbox{Im}(\gamma), (r_{ij})_{i\in\{0,1,2,3,4\}, j\in \{0,...,d\} })$ in $\overline{M}_{0,5d}$, and the coefficients $(c_i)_i$ give the $(\CC^*)^4$--fibration. The map $\gamma$ is determined by $(r_{ij})_{i,j}$ only up to quotient by the action of $(S_d)^{n+1}$, and the curve $\mbox{Im}(\gamma)$ is determined by $\gamma$ only after quotient by $\PP GL_2.$
\bp
For members of the Dwork pencil, the union $\bigcup_{\psi}\overline{M}_{0,0}(X_\psi,d,t)$ is a 625-1 cover over its image in $\overline{M}_{0,5d}$. 
\ep 
\bpr
The argument here was first mentioned for degree 1 maps in \cite{ak2}. Consider 
 the rational map on $\PP^4$ 
\bea f((x_0:...:x_4))=\frac{x_0^5+x_1^5+x_2^5+x_3^5+x_4^5}{x_0x_1x_2x_3x_4} 
\eea
The condition that $\mbox{Im}(\gamma) \subset X_\psi$ for some $\psi$ is equivalent to $ \gamma^*f$ being constant hence having the same zeroes and poles, with multiplicities.  In our case, for distinct roots $r_{ij}$, the condition reduces to \bea \sum_{i=0}^4 c_i^5 (r-r_{i1})^5\cdots (r-r_{id})^5 =0,  \quad \quad  \forall   r\in \{r_{ij} \}_{i\in\{0,1,2,3,4\}, j\in \{1,...,d\} }. \eea 
 Thus under suitable conditions on $r_{ij}$-s, the coefficients $(c_i^5)_i$  can be retrieved by solving a linear system of equations. 
 \epr
 In principle this could reduce the dimension of the ambient space for our space of degree $d$ curves in $X_\psi$ to $(5d-3)$. Even when $d=2,$ equations in 7 variables may be less than tractable.  Moreover,  this is only after taking quotients by  $\PP GL_2$ and $(S_d)^{n+1},$ and for degree $d>1$ the two actions do not seem to interact nicely. 
  
 Our approach to the space of conics will be different. We take from above the idea of working with curves with marked points, but instead of the marked points being traced by the hyperplanes of coordinates, we will focus on singularities of the associated plane sections (namely, intersection points of the conic with the associated cubic). Working in the spaces of planes with 2 such points allows us to reduce to the study of bitangent lines, and in the case of the Fermat quintic this reduces to a finite cover of $\PP^4.$ We are thus able to characterise the space of conics by four equations. These are rather involved and we have not yet started a systematic study exploring their geometric implications - we show here how they can be applied to identify specific families of conics in the Fermat quintic. 
 
 Finally, there are a number of reasons why one would be interested in the space of conics for the Fermat quitic and the Dwork pencil. With Clemens' conjecture in mind, it would be interesting to know how the higher degree curves interact with the loci of curves in degrees 1 and 2.  A distant but interesting goal concerns possible connections of these spaces to mirror symmetry of the  mirror quintic constructed from the Dwork pencil.
 %In \cite{cr} Herb Clemens and Ziv Ran proved

\section{Examples of Conics in Fermat quintic and the Dwork Pencil}

The group of permutations $S_5$ acting on $\PP^4$ leaves each Dwork quintic threefold $X_\psi$ invariant. We have found continuous families of conics in these threefolds by looking at planes kept invariant by certain subgroups of $S_5$. 

The first two examples exploit $S_3$ and $\ZZ_2\times \ZZ_2$--symmetries. These families were studied by Songyun Xu in \cite{xu}.

\begin{example} $S_3$-symmetry: \label{S3 conics}
Let $S\subset \PP^2\times \PP^1$
%$\rho: S\to \PP^1$ 
denote the 1-parameter family of quintic curves $S_\psi\hookrightarrow \PP^2$ given by
$$-5\psi^2ab^2c^2-5\psi a^3bc+a^5+b^5+c^5=0.$$
Let $P\subset \PP^4\times \PP^2$ denote the family of planes $P_{[a:b:c]}\subset \PP^4$ given by
    $$\begin{array}{cc} c(x_0+x_1+x_2)=ax_4,  & cx_3=bx_4\end{array}.$$
For each $([a:b:c], \psi)\in S$, the plane  $P_{[a:b:c]}$ intersects the quintic $X_\psi$ in the union of a conic and a cubic of equations: 
\bea && \mbox{ conic: } \quad x_0^2+x_1^2+x_2^2+x_0x_1+x_0x_2+x_1x_2-\psi x_3x_4 =0 \\  && \mbox{ cubic: } \quad (x_0+x_1)(x_0+x_2)(x_1+x_2)+\psi (x_0+x_1+x_2) x_3x_4=0. \eea 
\end{example}
Indeed, if we denote by $e_1, e_2, e_3$ the elementary symmetric polynomials in $x_0, x_1, x_2$ and set 
$$G_\psi :=x_3^5+x_4^5+e_1^5-5\psi x_3x_4e_1^3-5\psi ^2x_3^2x_4^2e_1$$ 
we get the following equation which allows us to factor $F_\psi$ under the conditions set up above:
 \bea
F_\psi -G_\psi =5(e_1^2-e_2-\psi x_3x_4)(e_3-e_1e_2-\psi e_1x_3x_4). 
\eea

We note that on the Fermat quintic, the cubic degenerates to a union of lines, while the conic is non-degenerate. The intersection points of the conic with the lines form the orbit under $S_3$  of the point $[a:-a: ia: ib: ic]$.

The following example was also found independently by Damiano Testa and Roger Heath-Brown by experimenting over finite fields. The plane sections here decompose into a union of one line and two conics. Johannes Walcher in \cite{walcher} extended this to families of conics on the Dwork pencil by considering the case when one $\ZZ_2$ subgroup swaps the two conics in each plane. 

\begin{example} \label{Z2 symmetries} {$\ZZ_2\times \ZZ_2$ symmetry for the Fermat quintic:}

For $i=1, 2$, let $u_i$ denote the elementary symmetric polynomials
in $x_0,x_1$, let $v_i$ denote the elementary symmetric polynomials
in $x_3,x_2$. Then
$$F=5u_1(u_2-\frac{u_1^2}{2})^2+5v_1(v_2-\frac{v_1^2}{2})^2+x_4^5-\frac{u_1^5}{4}-\frac{v_1^5}{4}.$$
It follows that, for the curve in $\PP^2$ given by the equation
$$ a^{10}+b^{10}-4b^5c^5=0,$$
the following are equations of conics in $X$:
 \bean \label{plane} && a^2(x_0+x_1)=b^2(x_2+x_3), \quad \quad bx_4=c(x_0+x_1) \quad \mbox{ and } \\ \quad && b(x_0^2+x_1^2)\pm i a(x_2^2+x_3^2)=0.\eean
 Under a change of variables $-4bc=d^2$ the parametrising curve can be brought to the nice form $a^{10}+b^{10}+d^{10}=0.$ The plane sections given by equations (\ref{plane}) also contain the special locus line identified in equation (\ref{isolated line }).
 %and the conics intersect this line at points ...

 \end{example}
 
In \cite{walcher}  Johannes  Walcher also extends this   method to look for conics invariant under $\ZZ_2$ generated by $(12)(34)\in S_5$, in planes of equations
\bea  a_1(x_1 + x_2) + a_2(x_3 + x_4) + x_5=0, \quad \quad 
(x_1 - x_2) + a_3(x_3 - x_4) =0 \eea
for suitable $a_i$-s. He found a family consisting entirely of reducible conics obtained by intersecting the van Geemen families of lines, and evidence of other solutions at the limit of computational capabilities.

\section{ Singular Plane Sections in the Fermat quintic}

We will now begin our systematic study of the space of conics in the Fermat quintic. We start by setting up the ambient space of singular plane sections. 

For the rest of the paper, $X$ will denote the Fermat quintic threefold. We use the notation $G(k,n)=G(\PP^{k-1}, \PP^{n-1})$ for Grassmannians.
%and $F(k_1,...,k_l,n)=F(\PP^{k_1-1 },...,\PP^{ -1 }, \PP^{n-1})$ for Flag varieties.

For a vector bundle $E$ 
%with locally free sheaf $\cE$, 
we will denote by $\PP(E)$
%$=\mbox{Proj}_X(\mbox{Sym}^*(\cE^\vee))$.
the projective bundle of lines in $E$.

% response to referee 1 
In line with the examples from the previous chapter, it is convenient to describe a plane in $\PP^4$ by the coefficients of two defining linear equations. 
To do so we consider the  canonical isomorphism $\CC^5 \cong (\CC^5)^\vee$ for a fixed choice of basis on $\CC^5$, which induces the isomorphism $G(3,5)\cong G(2,5)$ sending $\Lambda \subset \CC^5$ to $(\CC^5/\Lambda)^\vee \subset (\CC^5)^\vee$. Hence a plane $\Lambda$ given by independent equations
\bea \sum_{i=0}^4 a_ix_i=0 \quad \quad \mbox{and} \quad \quad  \sum_{i=0}^4 b_ix_i=0  \quad \quad \mbox{in } \quad \quad \PP^4 \eea 
will be described by 
\bea  \left[ \begin{array}{ccccc} a_0 & a_1 & a_2 & a_3 & a_4 \\ b_0 & b_1 & b_2 & a_b & a_b \end{array} \right] \in G(2,5)\eea
with Pl\"{u}cker coordinates $[a_ib_j-a_jb_i]_{0\leq i< j\leq 4}$ in $\PP(\bigwedge^2(\CC^5)^\vee)\cong \PP(\bigwedge^3\CC^5)=\PP^{9}$.

For each $\Lambda \in G(3,5)$, let $X_\Lambda$ denote the quintic plane curve  section $X_\Lambda = X \bigcap \Lambda.$

\begin{notation}
  We denote by $SP_k(X)$ the subset of  $G(3,5)$ made of planes for which $X_\Lambda$ has at least $k$ distinct singular points.
\end{notation}
%may upgrade to $SP_k(X)$
%\bea SP_k(X) = \{ \Lambda \in G(3,5); \quad l(X_\Lambda^s) \geq k\}  \eea  where for any variety $Y$, we let $Y^s$ be the locus of singular points of $Y$, and $l(Z)$ denotes the length of a 0-dimensional scheme $Z$.  \br Wait, do I need to localise at a point in order to define length? What is the terminology for the multiplicity of a point inside the singular locus? \er

$SP_k(X)$ comes with 
%a universal family 
an incidence variety $U_k\subset G(3,5)\times X^{k}$ given by
\bea U_k(X) = \left\{ (\Lambda, P_1,..., P_k) \in SP_k(X)\times X^{k}; \quad P_i \in X_\Lambda^s, \quad P_i \not= P_j \quad \forall i \not= j \right\}.\eea 
%\bea U_k(X) = \left\{ (\Lambda, P_1,..., P_k) \in SP_k(X)\times X^{k}; \quad P_i \in X_\Lambda^s, \quad \#\{i; P_i=p\}\leq l\left( (X_\Lambda^s)_P\right) \quad \forall P\in X_\Lambda^s \right\}\eea 
The  projection map $q_k: U_k(X) \to SP_k(X)$ forgets the marked points, and there is a natural group action of the permutation group $S_k$ on $U_k(X)$.
%which permutes the $k$ projection maps $\pi_i: U_k(X) \to X$.
%for $i\in \{1,...,k\}$.
  
\bp \label{U1} For $k=1$, the incidence variety above satisfies $U_1(X) \cong \PP(\Omega^1_X)$ where $\Omega^1_X$ is the cotangent sheaf of $X$.
As such, $U_1(X)$ can be thought of as a divisor in $\PP(E)$   for $E:=\cO_X(-1)^{\oplus 5}/ \cO_X(-5)$. 

Consider the rational map given in Pl\"{u}cker coordinates by
 \bea  X\times \PP^4 & \dra & G(2,5)  \\
       (P, y) & \dra & [p_i^4y_j-p_j^4y_i]_{0\leq i<j\leq 4} %[P^4\wedge Y]
       \eea
       where $P=[p_0:p_1:p_2:p_3:p_4]$ and $y=[y_0:y_1:y_2:y_3:y_4]$. 
       
       Let $q': \PP(\Omega^1_{\PP^4}|_X) \dra G(2,5)$ be the restriction of this map to %$\PP(\Omega^1_{\PP^4}\right)|_X)=V(\sum_{i=0}^4p_iy_i)\subset  X\times \PP^4.$
       \bea \PP(\Omega^1_{\PP^4}|_X)=\{(P,y)\in X\times \PP^4; \quad \sum_{i=0}^4p_iy_i=0\}. \eea
       %$\cong F(\PP^0,\PP^3,\PP^4)|_X  $
       
        Then the  map $q_1: U_1(X) \longrightarrow  SP_1(X) \subset G(3,5)\cong G(2,5)$ can be described as the vertical right hand side arrow in  the commutative diagram
        \bea \xymatrix{  \PP(\Omega^1_{\PP^4}|_X)  \ar@{-->}[r]^{\quad q'} \ar@{-->}[d]_{p'} &  G(2,5).  \\
\PP(\Omega^1_X) \ar[r]^\cong & U_1(X) \ar[u].  } \eea 
The  vertical arrow $p'$ is the dominant rational map
induced by the surjective homomorphism $\Omega^1_{\PP^4}|_X \to \Omega^1_X$. Both $p'$ and $q'$ have singular locus  $X'=\{(P,P^4)\in X\times \PP^4\} \cong X$ where $P^4:=[p_0^4:p_1^4:p_2^4:p_3^4:p_4^4]$.
\ep
   
  \bpr Indeed, a point $P$ is singular for a plane section $X_\Lambda$ if and only if $\Lambda \subset \ol{\cT}_{X,P}$ for the projective tangent space $\ol{\cT}_{X,P}$ of $X$ at $P$. Hence  such planes $\Lambda$ are parametrised by $\PP(\Omega^1_X)$.
  
  The cotangent sequence for $X$ and the restriction of 
  the Euler sequence for $\PP^4$ give the commutative diagram of vector bundles on $X$ with exact rows and columns:
  
 \bea \xymatrix{ & &0 \ar[d] & 0 \ar[d] & \\ 
0\ar[r] & \cO_X(-5)   \ar[r] \ar@{=}[d]&  \Omega^1_{\PP^4}|_X \ar[r]\ar[d] & \Omega^1_X \ar[r]\ar[d] & 0  \\
   0 \ar[r] & \cO_X(-5)    \ar[r]^{\nabla_F}  &\cO_X(-1)^{\oplus 5} \ar[r]\ar[d]  & E \ar[r]\ar[d] & 0 \\
 & & \cO_X  \ar@{=}[r]\ar[d]  & \cO_X  \ar[d] & \\
 & & 0 & 0 &
 }\eea
 where $\nabla_{F}(P)=(p_0^4, p_1^4, p_2^4, p_3^4, p_4^4)$ and $E:=\cO_X(-1)^{\oplus 5}/ \cO_X(-5).$  
 Hence $\PP(\Omega^1_X)$ is the divisor in $\PP(E)$  defined as the zero locus of the global section  $\sigma\in H^0(\cO_{\PP(E)}(1))$ given by
  \bea \cO_{\PP(E)} \longrightarrow p^*E \otimes  \cO_{\PP(E)}(1) \longrightarrow \cO_{\PP(E)}(1). \eea
  where $p:\PP(E) \to X$ is the projection map. Similarly, based on the middle vertical sequence we can describe $\PP(\Omega^1_{\PP^4}|_X)$ as the divisor in $\PP(\cO_X(-1)^{\oplus 5})=X\times \PP^4$ given concretely as the zero locus of $\sum_{i=0}^4p_iy_i.$
  
  The horizontal surjections in the commutative diagram above induce a dominant rational map $X\times \PP^4 \dra \PP(E)$
  and its restriction $\PP(\Omega^1_{\PP^4}|_X) \dra \PP(\Omega^1_X)$, with singular locus $\PP(\cO_X(-5))=X'.$
  For a pair $(P,y)\in (X\times \PP^4)\setminus X',$ the vectors $\nabla_{F}(P)$ and $(y_0,y_1,y_2,y_3,y_4)$ provide the coefficients for two independent linear equations defining a plane $\Lambda$ in $\PP^4.$ 
  Thus the rational map 
    $ X\times \PP^4 \dra G(2,5)$ given in  Proposition \ref{U1} is well defined on $(X\times \PP^4)\setminus X'$ and constant on the fibres of $(X\times \PP^4)\setminus X' \to \PP(E)$, hence descends to a morphism $\PP(E) \to G(2,5).$ By restriction, the rational map $q': \PP(\Omega^1_{\PP^4}|_X) \dra G(2,5)$ is well defined on $\PP(\Omega^1_{\PP^4}|_X)\setminus X'$ and constant on the fibres of $\PP(\Omega^1_{\PP^4}|_X)\setminus X' \to \PP(\Omega^1_X)$, descending to a morphism $\PP(\Omega^1_X) \to G(2,5).$ Composition with the isomorphisms $ U_1(X)\cong \PP(\Omega^1_X)$ and $G(2,5)\cong G(3,5)$ yields exactly $q_1: U_1(X) \longrightarrow  SP_1(X) \subset G(3,5).$ Indeed, with the notations above we have $\Lambda \subset  \ol{\cT}_{X,P}$ whenever $(P,y)\in  \PP(\Omega^1_{\PP^4}|_X)\setminus X'$.
  \epr

\begin{notation}
  We denote by $\Delta_{2,2}(X)$ the variety of bitangent lines with marked tangency point: 
  \bea  \Delta_{2,2}(X)=\left\{ (l, P, Q)\in G(2,5)\times X\times X; \quad P \not= Q, \quad l\cdot X \geq 2P+ 2Q \right\} \eea
  which parametrizes lines tangent to $X$ at two (distinct) marked points  (compare with \cite{cr}).
  
\end{notation}

  For every space $Y$, we define $\Delta_Y$ as the diagonal $\Delta_Y:=\mbox{Im}(Y \hookrightarrow Y\times Y)$. 
  
In the following discussion we will identify $\ZZ^4_4$ with the quotient of the group homomorphism $i: \ZZ_4 \to \ZZ_4^5$ given by $i(k_0)=(k_0,k_0,k_0,k_0,k_0)$, and as such we will write elements of $\ZZ_4^4$ as equivalence classes $[k]:=[{k_i}]_{0\leq i \leq 4}$ with $k_i\in\ZZ_4.$

Consider the action of $\ZZ_4^4$ on $\PP^4\times \PP^4$    \bean\label{Z4 locus }\sigma : \ZZ_4^4 \times (\PP^4\times \PP^4) & \longrightarrow & \PP^4\times \PP^4 \\
\left([k],  \left([p_i]_i,  [q_i]_i \right)\right) & \longrightarrow & \left( [p_i]_i,  [\sqrt{-1}^{k_i}q_i]_i  \right),
  \eean
  and let $\sigma_{[k]}=\sigma([k], \_\_)$ for each $[k]=[{k_i}]_{0\leq i \leq 4}$ as above.
 Let $X_{[k]}:=\sigma_{[k]}(\Delta_{\PP^4})\bigcap (X\times X)$ and $X_{\sigma}:=\ZZ_4^4\cdot(\Delta_{\PP^4})\bigcap (X\times X)=\bigcup_{[k]\in\ZZ_4^4}X_{[k]}$.

\bt \label{U2} We have a commutative diagram 
\bea \xymatrix{  {U}_2(X)  \ar[d]^{f^\circ} \ar[r]^{\cong \quad \quad \quad \quad \quad \quad} &\left(\PP(\Omega^1_X)\times_{G(3,5)}\PP(\Omega^1_X)\right) \setminus \Delta_{\PP(\Omega^1_X)} \ar[d] \ar@{^{(}->}[r] & \PP(\Omega^1_X)\times_{G(3,5)}\PP(\Omega^1_X) \ar[d]^p \\ \Delta_{2,2}(X) \ar@{^{(}->}[r] & (X \times X)\setminus \Delta_X \ar@{^{(}->}[r] & X\times X  } \eea 
such that the vertical map $f^\circ$ extends to a surjective morphism $f:\overline{U}_2(X)\ra \overline{\Delta}_{2,2}(X)$, where 
 $\overline{\Delta}_{2,2}(X)$ is the closure of $\Delta_{2,2}(X)$ in $X \times X$ 
 and  $\overline{U}_2(X)$ is the closure of $U_2(X)$ in $\PP(\Omega^1_X)\times_{G(3,5)}\PP(\Omega^1_X).$ 
% Both $\overline{\Delta}_{2,2}(X)$ and $\overline{U}_2(X)$ are 4-dimensional.

For every $[k]\in \ZZ_4^4$ non-trivial, $X_{[k]}$ is a 2-dimensional 
%complete intersection 
subvariety of $\overline{\Delta}_{2,2}(X)$ and $f$ has $\PP^1$--fibres over all points of  $X_{[k]}\setminus X_{[0]}$. Outside $X_\sigma=\bigcup_{[k]\in\ZZ_4^4}X_{[k]}$, f is an isomorphism.
%\br We probably need at least $ \mbox{Bl}_{\Delta_X}(X\times X)$ in order to say something about what happens over $\Delta_X$? In this case we should be able to write the exceptional locus in coordinates... do we need to blow up all $X_\sigma$ or component by component? Some of these are not smooth. careful on the interaction between the blow-up locus and the diagonal.\er
 
 \et
 
  \bpr  The two projections $U_2(X)\to U_1(X)$ together with the isomorphism in Proposition \ref{U1} give an open embedding $U_2(X) \ra \PP(\Omega^1_X)\times_{G(3,5)}\PP(\Omega^1_X)$, whose  image is exactly $\left(\PP(\Omega^1_X)\times_{G(3,5)}\PP(\Omega^1_X)\right) \setminus \Delta_{\PP(\Omega^1_X)}$.
  %The horizontal arrows in the diagram follow directly from the definitions. 
  
  Given a plane $\Lambda \in G(3,5)$, two distinct points $P$ and $Q$ are singular in the plane section $X_\Lambda$ iff $PQ\subset \Lambda \subseteq \ol{\cT}_{X,P} \bigcap  \ol{\cT}_{X,Q}.$ Hence $(PQ, P, Q)\in \Delta_{2,2}(X)$. This defines $f^\circ.$

For generic $(PQ, P, Q)\in \Delta_{2,2}(X)$ we have a unique plane $\Lambda$ 
  %such that $(\Lambda, P, Q)\in U_2(X)$,
  as above, namely $\Lambda = \ol{\cT}_{X,P} \bigcap  \ol{\cT}_{X,Q}$. The special points $(P,Q)\in X\times X$ where $\ol{\cT}_{X,P} =  \ol{\cT}_{X,Q}$ are those for which  \bea [\nabla_F(P)]:=[p_0^4:p_1^4:p_2^4:p_3^4:p_4^4]=[q_0^4:q_1^4:q_2^4:q_3^4:q_4^4]=:[\nabla_F(Q)], \eea namely those pairs in $X\times X$ which are also in $\ZZ_4^4\cdot(\Delta_{\PP^4}).$ In other words, the restriction of $f^\circ$ on the preimage of $\Delta_{2,2}(X) \setminus X_\sigma$  gives a bijective map to $\Delta_{2,2}(X) \setminus X_\sigma$.
  %Here $i\in\{0,1,2,3,4\}$ and $k_i\in\{0,1,2,3\}$.
  
     $\Delta_{2,2}(X)$ is a closed subset of $(X \times X)\setminus \Delta_X$  given in $(\PP^4\times \PP^4) \setminus \Delta_{\PP^4}$ by equations \bean \label{double tangents} \sum_{i=0}^4 q_i^hp_i^{5-h} =0 \mbox{ for }   h\in \{0,1,4,5\}\eean 
for   distinct $P=[p_0:p_1:p_2:p_3:p_4]$ and $ Q=[q_0:q_1:q_2:q_3:q_4]\in X$, hence its components will have dimension at least 4.

%\textbf{In multi-homogeneous coordinates, the Jacobian matrix\bea \left( \begin{array}{cccccccccc} 5p_0^4 &  5p_1^4 &  5p_2^4 &  5p_3^4 &  5p_4^4 & 0& 0& 0& 0& 0 \\ 4p_0^3q_0 &  4p_1^3q_1 &  4p_2^3q_2 &  4p_3^3q_3 &  4p_4^3q_4 & p_0^4 & p_1^4 &  p_2^4 &  p_3^4 &  p_4^4 \\  q_0^4 & q_1^4 &  q_2^4 &  q_3^4 &  q_4^4 &4p_0q_0^3 &  4p_1q_1^3 &  4p_2q_2^3 &  4p_3q_3^3 &  4p_4q_4^3  \\0& 0& 0& 0& 0 & 5q_0^4 & 5q_1^4 &  5q_2^4 &  5q_3^4 &  5q_4^4  \end{array} \right) \eea has rank 4 on all points in $(\PP^4\times \PP^4) \setminus \ZZ_4^4\cdot\Delta_{\PP^4}$}
  
  In the case when $Q=[k]\cdot P$ for some $[k]=[\sqrt{-1}^{k_i}]_i\in \ZZ_4^4$ non-trivial, the equations (\ref{double tangents}) for $h=0,4$ and for $h=1,5$ coincide, hence $X_{[k]}\hookrightarrow \overline{\Delta}_{2,2}(X)$ and $X_{[k]}$ is isomorphic to the surface $V\left(\sum_{i=0}^4 p_i^5, \sum_{i=0}^4 p_i^5\sqrt{-1}^{k_i}\right)$ in $X$, while $X_{[0]}\cong X$. Thus all components of $X_\sigma$ have dimension smaller than $4$ and hence $f^\circ$ is a dominant map. 
  
  With the identifications above, $f^\circ$ is just the restriction to $\left(\PP(\Omega^1_X)\times_{G(3,5)}\PP(\Omega^1_X)\right) \setminus \Delta_{\PP(\Omega^1_X)}$ of the projection map $p: \PP(\Omega^1_X)\times_{G(3,5)}\PP(\Omega^1_X) \to X\times  X$ which remembers the two marked points coming from the two copies of $\PP(\Omega^1_X)\cong U_1(X)$. As such, $f^\circ$ extends to $f:\overline{U}_2(X)\ra \overline{\Delta}_{2,2}(X)$, which is just the restriction of $p$ on $\overline{\left(\PP(\Omega^1_X)\times_{G(3,5)}\PP(\Omega^1_X)\right) \setminus \Delta_{\PP(\Omega^1_X)}}.$ Since  $\ol{U}_2(X)$ and $\ol{\Delta}_{2,2}(X)$ are compact and $f^\circ$ is dominant, we get that $f$ is a surjective morphism. In fact, we have seen that $f$ is invertible outside $X_\sigma=\bigcup_{[k]\in\ZZ_4^4}X_{[k]}$.
  
   For $(P,Q)\in X_{[k]}$ with $[k]\in\ZZ_4^4$ non-trivial, the fibre $f^{-1}((P,Q))$ is the pencil $\PP^1$ parametrising planes $\Lambda$ with $PQ\subset \Lambda \subset \ol{\cT}_{X,P}= \ol{\cT}_{X,Q}.$
  \epr

   The surjective map $q_2: U_2(X) \to SP_2(X)$ can be extended to $q_2: \overline{U}_2(X) \to \overline{SP}_2(X)$ by restricting the projection $q:\PP(\Omega^1_X)\times_{G(3,5)}\PP(\Omega^1_X) \ra G(3,5)$.
 The same map can also be described as the composition \bea \overline{U}_2(X) \ra X\times X \dra G(3,5),\eea where the right side map is the restriction of the rational map $g: \PP^4\times \PP^4 \dra G(3,5)$ which via the isomorphism $G(3,5)\cong G(2,5)$ is defined in Pl\"{u}cker coordinates by 
   \bea\xymatrix{\left([p_0:p_1:p_2:p_3:p_4], [q_0:q_1:q_2:q_3:q_4] \right) \ar[r]^{ \quad \quad \quad \quad g}  & [p_i^4q_j^4-p_j^4q_i^4]_{0\leq i<j\leq 4}.}\eea
     We note that this map can be resolved by blowing up along the locus $\ZZ_4^4\cdot\Delta_{\PP^4}$ discussed above. 
  
\subsection{Coordinates for $\overline{\Delta}_{2,2}(X)$ } 
\begin{notation} \label{u coords}
  
To set up our new coordinates, consider the rational map $\phi: \overline{\Delta}_{2,2}(X) \dashrightarrow \PP^4$ obtained as a composition $\phi=u\circ\pi$ between the projection $\pi: \overline{\Delta}_{2,2}(X) \rightarrow \PP^4\times\PP^4$ and the rational map 
   \bea  u: \PP^4\times\PP^4 \dashrightarrow \PP^4, \quad \quad  
  %over $U\times \PP^4$ 
   u(P,Q)= U:=[u_0:u_1:u_2:u_3:u_4]\eea where $u_i=q_i/p_i$ for all $i.$
  \end{notation}
 %   \br Resolve phi: is this ok? The strict transform of $\Delta_{2.2} $ should be taken in the blow-up of $\ol{\c{LM}_{0,7}} \times \ol{\c{LM}_{0,7}}$ along all $E_i \times E_i$ where $E_i$ are the boundary divisors in $\ol{\c{LM}}_{0,7}$ \er 
    
    Consider the involution  $\tau : \overline{\Delta}_{2,2}(X) \rightarrow \overline{\Delta}_{2,2}(X)$  which swaps $P$ and $Q$. There is a corresponding map $\iota: \PP^4 \dashrightarrow \PP^4$ such that the following diagram  is commutative  
  \bea \xymatrix{  \overline{\Delta}_{2,2}(X)  \ar@{-->}[d]_\phi \ar[r]^\tau & \overline{\Delta}_{2,2}(X) \ar@{-->}[d]_\phi  \\
 \PP^4 \ar@{-->}[r]^\iota  & \PP^4  } \eea
%used manual at https://personal.math.ubc.ca/~cautis/tools/xypic.pdf
%where for $U=[u_0:u_1:u_2:u_3:u_4]$
The map $\iota: \PP^4 \dashrightarrow \PP^4$ is the Cremona transformation given on an open set by $u_i \ra 1/u_i.$

We will now prepare some notations for more concrete coordinate calculations:

\begin{notation} \label{symm notations}
  The discriminant and elementary symmetric polynomials
  \bea \delta :=\prod_{j>k} (u_j-u_k), \quad  e_0=1, \quad   e_1=\sum_{j=0}^4u_j,   \quad  e_2=\sum_{4\geq j>k\geq 0}u_ju_k, \quad ... \quad, e_5=\prod_{j=0}^4u_j\eea
  will be useful in our computations. We will also work with the $S_4$--symmetric polynomials
\bea e_0(i)=1, \quad \quad  e_1(i)=\sum_{j\in \{0,...,4\}\setminus\{i\}}u_j,   \quad \quad e_2(i)=\sum_{j,k\in \{0,...,4\}\setminus\{i\}\mbox{; }j>k}u_ju_k, \quad \mbox{ etc. } \eea 
and the partial discriminants
 \bea  d_i :=\prod_{j,k\in \{0,...,4\}\setminus\{i\}\mbox{; }j>k} (u_j-u_k).\eea
   
\end{notation}

  Let $\varphi: X \to H\cong \PP^3$ denote the quotient map for the action of $(\ZZ_5)^4$ on $X$,
  \bea \varphi([x_0:x_1:x_2:x_3:x_4]) = [x_0^5:x_1^5:x_2^5:x_3^5:x_4^5],\eea which maps $X$ to the hyperplane $H=V(\sum_{i=0}^4y_i)$ in $\PP^4=\mbox{Proj}\CC[y_0,y_1,y_2,y_3,y_4].$ We note that $\varphi$ is ramified exactly along the coordinate hyperplanes $H_i=V(x_i)$.

  \bt \label{P4 parametrisation} The map $\phi=u\circ\pi: \overline{\Delta}_{2,2}(X) \dashrightarrow \PP^4$  is generically $625:1$, unramified over its regular points.  There is a commutative diagram   
  \bea \xymatrix{  \overline{\Delta}_{2,2}(X) \ar@{-->}[d]_\phi \ar[r]^{\quad \pi_1} & X \ar[d]^\varphi \\
 \PP^4\ar@{-->}[r]^{m}  & H  } \eea
%used manual at https://personal.math.ubc.ca/~cautis/tools/xypic.pdf
where for $U=[u_0:u_1:u_2:u_3:u_4]$, we have $m(U)=[M_i(U)]_{i\in\{0,...,4\}}$ given by
\bea M_i(U) := (-1)^id_in_i, \quad \mbox{ where } \quad n_i := e_2(i)^2-e_1(i)e_3(i).\eea
Along the horizontal rows have $\pi_1\circ \tau = \pi_2$ and $m\circ\iota(U)=[M_i(U)u_i^5]_{i\in\{0,...,4\}}$.

The base locus of $m$ is  made of 25 planes $V(u_i-u_j,u_k-u_l)$ for $i\not=j$,  $k\not=l$ and $\{ i, j\} \not=\{ k, l\}$ (of which 10 correspond to the case $\left|\{i, j\} \bigcap \{k, l\}\right|=1$ and 15 to the case $\{i, j\} \bigcap \{k, l\}=\emptyset$), 10 irreducible quartic surfaces  $V(u_i-u_j,n_i)$ for $i>j$, and a degree 6 irreducible component $V(e_2,e_3)$. 
  \et
  
  \bpr
    Let $V\cong (\CC^*)^4$ and $W $ denote the main open sets in $ \PP^4$ and $G(3,5)$ respectively, obtained by removing the coordinate hyperplanes (for the Pl\"{u}cker coordinates in the case of $G(3,5)$).
  
  In projective coordinates $P=[p_0:p_1:p_2:p_3:p_4]$ and $Q=[q_0:q_1:q_2:q_3:q_4]$ and  \bea PQ=\{tP+sQ=[tp_0+sq_0:tp_1+sq_1:tp_2+sq_2:tp_3+sq_3:tp_4+sq_4]; \quad [t:s]\in \PP^1  \},\eea the condition $(PQ, P, Q)\in \Delta_{2,2}(X)$ translates into the four equations (\ref{double tangents}).
 With the new coordinates introduced in Notation \ref{u coords}, equations (\ref{double tangents}) are written as
    \bean \label{eq Mi} \left( \begin{array}{ccccc} 1&1&1&1&1 \\
                          u_0 & u_1 & u_2 & u_3 & u_4 \\
                          u_0^4 & u_1^4 & u_2^4 & u_3^4 & u_4^4 \\
                          u_0^5 & u_1^5 & u_2^5 & u_3^5 & u_4^5 
 \end{array} \right)  \left( \begin{array}{c} p_0^5 \\p_1^5\\p_2^5 \\ p_3^5 \\ p_4^5 \end{array} \right) =
  \left( \begin{array}{c} 0 \\ 0 \\ 0 \\ 0 \\ 0 \end{array} \right) \eean
  which admits the solution $ p_i^5=M_i(U)$, for $M_i(U)$ defined as in the statement of the Theorem. This defines the map $m$. The base locus of $m$ is  $V(\left\{M_i(U)\right\}_{i\in\{0,...,4\}})$, which consists of 25 planes $V(u_i-u_j,u_k-u_l)$ for $i< j\leq k< l$, 10 irreducible quartic surfaces  $V(u_i-u_j,n_i)$ for $i>j$, and a degree 6 irreducible component $V(e_2,e_3)$. 
  
  Indeed,  $(u_i-u_j)$ is a factor of three sections $M_k(U)\in H^0(\cO_{\PP^4}(10))$, and also of $M_i(U)-M_j(U)$. This explains the first two types of components. 
  
  The last component is in fact  $V(\left\{ n_i \right\}_{i\in\{0,...,4\}})=V(e_2,e_3)$. Indeed, on the one hand
  \bea e_2=e_2(i)+u_i e_1(i) \quad \mbox{ and } \quad e_3=u_i e_2(i)+e_3(i)\eea 
   so $V\left(e_2, e_3\right) \subseteq V\left( n_i\right)$. In addition, within the ideal $\left( n_i  \right)_{i\in\{0,...,4\}}$ we have the identities
  \bea n_i(u_j-u_k)+n_j(u_k-u_i)+n_k(u_i-u_j)=(u_j-u_k)(u_k-u_i)(u_i-u_j)e_2 \eea for all triples $(i,j,k)$ while the following symmetric sums satisfy
  \bea \sum_{i=0}^4n_i=3e_2^2-4e_1e_3, \quad  \sum_{i=0}^4u_in_i=e_2e_3, \quad  \sum_{i=0}^4u_i^2n_i=3e_3^2-4e_4e_2,\eea
hence $V\left(e_2, e_3\right) = V(\left\{ n_i \right\}_{i\in\{0,...,4\}})$.

The primary decomposition of the ideal $(\left\{M_i(U)\right)_{i\in\{0,...,4\}})$, checked with the help of Macaulay2, shows that all the components of the base locus of $m$ are accounted for by the list stated in the Theorem.
There are no embedded components, hence the map $m$ can be resolved by successive blow-ups of $\PP^4$ along the loci made of these components and their intersections, starting with the smallest dimension one. 
Blow-ups of $V(\left\{M_i(U)\right\}_{i\in\{0,...,4\}})$ are part of the construction of $\ol{\cM}_{0,7}$ as discussed below. The other components of the singular locus above need to be considered separately.

We note that for the involution $\iota$  \bea e_k(\iota(U)) e_5=e_{5-k}(U), 
%\quad e_k(i)(\iota(U)) e_5=e_{4-k}(i)(U)u_i, 
\quad d_i(\iota(U)) e_5^3=d_i(U) u_i^3,  \quad n_i(\iota(U)) e_5^2=n_i(U) u_i^2  \eea
hence $M_i(\iota(U))e_5^5= M_i(U) u_i^5$ which is indeed the descent of the relation $\pi_1\circ\tau=\pi_2.$
 \epr
 % see where in the future this will occur

%  \br We will need to look separately at the cases of planes in the boundary of the Grassmannian? These should be planes intersecting "lines of coordinate" $V(x_i,x_j,x_k)$. \er  \br I need to come back and say a few words about more singularities and their relations to the map m. This might be handy when discussing special cases. \er 

\subsection{Resolving the maps $\phi$, $\tau$ and $m$. } \label{resolution}

Let $H_i=V(x_i)$ denote the coordinate hyperplanes in $\PP^4$. 
To resolve both $\iota$ and $\varphi$ simultaneously, it is enough to blow-up $\PP^4$ along $\bigcap_{i\in J}H_i$ for all $J\subset \{0,1,2,3,4\}$ with $\#J>1,$ done successively from smallest to largest strata. The resulting space is a Losev-Manin moduli space of rational curves with marked points. To partially resolve $m$, we can consider $\ol{\cM}_{0,7}$, obtained by successive blow-ups of $\PP^4$ along $\bigcap_{\alpha\in J}H_\alpha$ for $J\subset \{0,1,2,3,4\}\bigcup \{(i,j); \quad 0 \leq i<j\leq 4\}$, where $H_{(i,j)}=V(x_i-x_j)$. Let $\widetilde{\Delta}_{2,2}(X)$  denote the strict transform of  $\overline{\Delta}_{2,2}(X)$ in  $\ol{\cM}_{0,7}\times \ol{\cM}_{0,7}$. We then obtain a lift of the earlier diagram to 
\bea \xymatrix{  \widetilde{\Delta}_{2,2}(X)  \ar[d]_{\widetilde\phi} \ar[r]^{\widetilde\tau} & \widetilde{\Delta}_{2,2}(X) \ar[d]_{\widetilde\phi}  \\
\ol{\cM}_{0,7}\ar[r]^{\widetilde\iota}  & \ol{\cM}_{0,7}.  } \eea

This also partially resolves the map $m$, though there remains the locus at the intersection of $n_i$-s which would need a separate blow-up.

\section{Parametrising Conics with high contact at  Singular points of Plane Sections of the Fermat quintic}

Consider the singular plane sections of the Fermat quintic parametrised by $\ol{SP}_2(X)$, and their incidence variety $\ol{U}_2(X)$. We fix an element $(\Lambda, P, Q)\in {U}_2(X)$. For the next steps we assume $(PQ, P, Q)\in \Delta_{2,2}(X)\setminus X_{\sigma}$, though in the last sections we will discuss how the methods developed here can be adapted to points in $ X_{\sigma}$.

We would like to find conditions for a conic $C\hookrightarrow \Lambda$ to be a component of $X_\Lambda$. We will again approach this problem gradually.

 We write the fixed points as $P=[p_0:p_1:p_2:p_3:p_4]$,  $Q=[q_0:q_1:q_2:q_3:q_4]$ and $\phi(PQ, P, Q)=U=[u_0:u_1:u_2:u_3:u_4]$.  (Even though it is convenient to use homogeneous coordinate notations, we fix affine representatives for the above).
  In the previous section we have seen how, under the condition $(PQ, P, Q)\in \Delta_{2,2}(X)$, the points $P$ and $Q$ can be recovered (though only up to the action of the group $\ZZ_5^4$) from the data $U=[u_i]_i$, via the formulae $[p_i^5]_{i\in\{0,...,4\}}:=m(U)=[M_i(U)]_{i\in\{0,...,4\}}$ and $q_i^5:=M_i(U)u_i^5$ given in Theorem \ref{P4 parametrisation}. Furthermore, if the section $X_\Lambda$ contains a conic, then in general we expect to find other singularities $S\in X_\Lambda^s$, so that $(\Lambda, P, S)\in {U}_2(X)$. Then similarly for $S=[s_i]_i$ we could consider new coordinates $v_i=s_i/p_i$, such that up to  the action of the group $\ZZ_5^4$, both $P$ and $S$ can be written in terms of $V=[v_i]_i$. The search for such $S$ motivates us to perform a change of coordinates as follows:
 
 \begin{notation}\label{nup}
 
 Consider the map 
 \bea \nu_P: \PP^4 \dra \PP^4, \quad \quad  \nu_P([v_i]_{i\in\{0,1,2,3,4\}}) = [p_iv_i]_{i\in\{0,1,2,3,4\}}. \eea
   
 \end{notation}
 Note that here the values $\{p_i\}_i$ are fixed. 
  For the next steps we will assume that $P\in V=\PP^4\setminus \bigcup_iH_i$ so that $\nu_P$ is an isomorphism. We will work on the domain of $\nu_P$, namely on $\PP^4=\mbox{Proj }\CC[v_0,v_1,v_2,v_3,v_4],$ and describe $X$, $\Lambda$ as well as $P$, $Q$ via the change of coordinates $\nu_P^{-1}([x_i]_i)=[x_i/p_i]_i$. The coordinates $x_i$ are pulled-back as $\nu_P^*(x_i)=p_iv_i$, and $P$ is identified with $\nu_P^{-1}(P)=[1]:=[1:1:1:1:1]$, and $Q$ with $\nu_P^{-1}(Q)=[U]=[u_0:u_1:u_2:u_3:u_4]$. 
  
  Assume $u_i\not=u_j$ for $i\not=j$. Let $F(x)=\sum_{i=0}^4x_i^5\in H^0(\cO_{\PP^4}(5)).$ We identify the Fermat quintic $X=V(F)$ with $\nu_P^{-1}(X)=V(\nu_P^*(F))$ where in the new coordinates  
  \bean \label{new quintic eqs} \nu_P^*(F)(v) = \sum_{i=0}^4p_i^5v_i^5 = \sum_{i=0}^4M_i(U)v_i^5.\eean 
  Similarly, for the plane $\Lambda$, defined in the proof of Theorem \ref{U2} by 
  \bea \Lambda=\ol{\cT}_{X,P} \bigcap  \ol{\cT}_{X,Q}=V\left(\sum_{i=0}^4p_i^4x_i, \sum_{i=0}^4q_i^4x_i\right),\eea
  we will now identify $\Lambda$ with $\nu_P^{-1}(\Lambda)$
  given by equations
  % given by $k=1,4$ in formulas (\ref{double tangents}), have become
   \bean \label{eqs for lambda} \nu_P^*(\sum_{i=0}^4p_i^4x_i) =\sum_{i=0}^4M_i(U)v_i = 0  \mbox{ and }   \nu_P^*(\sum_{i=0}^4q_i^4x_i) =\sum_{i=0}^4M_i(U)u_i^4v_i = 0, \eean
%  \bean \label{new plane eqs1} & & \sum_{i=0}^4M_i(U)v_i = 0 \quad \quad \mbox{ and }\\ \label{new plane eqs2} & & \sum_{i=0}^4M_i(U)u_i^4v_i = 0, \eean 
  Hence the new coordinates allow us to write all important equations based on the fixed data $[U]=[u_i]_i$, in which terms we already know two fixed singular points $P, Q\in X_\Lambda$. 
  
  We also note that the involution $\tau$ swapping $P$ and $Q$ generates the following changes: 
  \bean \label{involution} U=[u_i]_i \ra \iota(U)=\left[\frac{1}{u_i}\right]_i, \quad    M_i(U) \ra M_i(\iota(U))=\frac{M_i(U)u_i^5}{e_5^5},  \quad v_i \ra \frac{v_i}{u_i}  \eean 
   for $ i\in\{0,1,2,3,4\}$.
   
We will now start looking for conditions that insure 
\bea  C\cdot X_{\Lambda}\geq 3P + 3Q \quad \mbox{ in }  \Lambda. \eea 
   
   \begin{notation}
       We now consider the conics in $\PP^4$ given by
  \bean \label{conic parametrisation} \gamma_{b}: \PP^1 \ra \PP^4, \quad \quad \quad \quad \gamma_{b}([t:s])= [t^2+b_i ts + u_ib_5^2s^2]_{i}, \eean 
   so that $\gamma_{b}([1:0])=P$ and $\gamma_{b}([0:1])=Q$, and depending on  $b=(b_0,b_1,b_2,b_3,b_4,b_5)\in \CC^6.$ We will denote $C_{[b]}=\mbox{Im}(\gamma_{b})$.
   \end{notation}
  
 The parameter $b$ was chosen so that under the action of the remaining automorphism group $\CC^*=\mbox{Spec } \CC[\lambda, \frac{1}{\lambda}]\subset \PP GL_2,$ 
  \bea \gamma_{b}([\lambda t:s])= [t^2+\frac{b_i}{\lambda} ts + u_i\left(\frac{b_5}{\lambda}\right)^2s^2]_{i}= \gamma_{b/\lambda}([t:s]), \eea 
  so that we have a $\PP^5$--family of conics $C_{[b]}$ with $[b]\in \PP^5$. Let $o=(0,0,0,0,0,1)$.

\bp Consider the $(\PP^5\setminus\{[o]\})$--family of conics passing through $P$ and $Q$ with the parametrisation (\ref{conic parametrisation}) above. Let $B:=[b_0:b_1:b_2:b_3:b_4]$. Then the condition
\bea  C_b \subset \Lambda \mbox{  and } C_b\cdot X_{\Lambda}\geq 3P + 3Q \mbox{ in } \Lambda \eea
is equivalent to  $B \in \Lambda \bigcap \Gamma_U \bigcap \Gamma_{\iota (U)}$, where $\Gamma_U$ and  $\Gamma_{\iota (U)}$ are the quadric surfaces in $\PP^4$ given by equations
\bean \label{tangent cone P} & & \sum_{i=0}^4M_i(U)v^2_i = 0  \quad\quad \quad \mbox{ for } \quad \quad \Gamma_U \quad \quad \mbox{ and } \\ & & \label{tangent cone Q} \sum_{i=0}^4M_i(U)u_i^3v^2_i  =  0  \quad \quad \mbox{ for } \quad \quad \Gamma_{\iota (U) }.\eean 
\ep

\bpr We have chosen the conic parametrisation in (\ref{conic parametrisation}) so that $B$ represents its tangent direction at $0=[0:1]$ and $\infty=[1:0].$ The incidence conditions at $P$ and $Q$ above are equivalent to $B$ being in the tangent cones to $X_\Lambda$ at $P$ and $Q$. The equations (\ref{tangent cone P}) and (\ref{tangent cone Q}) above define the two tangent cones. In coordinates, we can write
% We evaluate $F$ at $\gamma_b:$
  \bea f_b(t)=\nu_P^*(F)(\gamma_{b}([t:1]))=\sum_{i=0}^4M_i(U)(t^2+b_i ts + u_ib_5^2s^2)^5. \eea Then the conditions for $C_b \subset \Lambda $  and $ C_b\cdot X_{\Lambda}\geq 3P + 3Q$ in $\Lambda$   can be calculated as
    \bea  \eval{\frac{\partial^k f_b}{\partial t^k}}_{t=0}=0 \quad \quad  \mbox{ for } k\in \{0, 1, 2, 8, 9, 10\}.\eea
Conditions for $(10-k)$ are equivalent to $\eval{\frac{\partial^k f_b}{\partial t^k}}_{t=\infty}=0$ and can be obtained from conditions for $k$ by extending the transformation $\iota$: 
\bea U=[u_i]_i \ra \iota(U)=\left[\frac{1}{u_i}\right]_i, \quad    M_i(U) \ra M_i(\iota(U))=\frac{M_i(U)u_i^5}{e_5^5},  \quad b_i \ra \frac{b_i}{u_i} \mbox{ for } i\in\{0,...,4\}. \eea 
%  Indeed, $\gamma_{b}([t:s])= [t^2+b_i ts + u_ib_5^2s^2]_{i}= [\frac{1}{u_i}t^2 +\frac{b_i}{u_i} ts + b_5^2s^2]_{i} =\gamma_{(1, (b_i/u_ib_5)_i)}([s:t/b_5])$.
We obtain the following set of conditions for $k\in\{0,1,2\}$ at $t=0$ and $t=\infty$:   
  \bea \begin{array}{llll}
    k=0: & \quad \quad \sum_{i=0}^4M_i(U) = 0 \quad & \mbox{ and } & \quad \quad  \sum_{i=0}^4M_i(U)u_i^5 = 0,\\
    k=1: & \quad \quad \sum_{i=0}^4M_i(U)b_i = 0  \quad \quad & \mbox{ and }&  \quad \quad \sum_{i=0}^4M_i(U)u_i^4b_i = 0,  \\   k=2: & \quad \quad \sum_{i=0}^4M_i(U)b^2_i  =  0 \quad \quad & \mbox{ and }&  \quad \quad \sum_{i=0}^4M_i(U)u_i^3b^2_i  =  0. \end{array}\eea
  The identities for $k=0$ correspond to $P,Q\in X$ and impose no conditions on the family of conics. The identities for $k=1$ are equivalent to $B:=[b_0:b_1:b_2:b_3:b_4]\in \Lambda$ and the cases $k=2$  to $B\in \Gamma_U \bigcap \Gamma_{\iota (U)}$, the intersection of the two tangent cones.
  \epr 
  
  Next we will choose a third point $R\in \Lambda\setminus \{PQ\}$ and thus set a system of coordinates $[x:y:z]$ for any point $V=[xP+yQ+zR]\in\Lambda$.

  \bl  In general, the conic $\Gamma_U \bigcap \Lambda$ is a union of two lines intersecting at $P$ while $\Gamma_{\iota (U)}\bigcap \Lambda$ is a union of two lines intersecting at $Q$. \el
  
  \bpr Indeed, for every $b\in \PP^5$ whose projection on $V(b_5)\cong \PP^4$ is $B\in \Gamma_U \bigcap \Lambda$, the indexed conic $C_b$ must satisfy
\bea  C_b\cdot X_{\Lambda}\geq 3P + 2Q \quad \mbox{ in }  \Lambda. \eea   
  Since $X_{\Lambda}$ is singular at both $P$ and $Q$, in general we expect that the tangent cone to $X_{\Lambda}$ at $P$ is a union of lines $l^1_P$ and $l^2_P$, and similarly for the tangent cone at $Q$. Then each conic $C_b$ must pass through $Q$ and be tangent to one of $l^1_P$ or $l^2_P$. For each $l^i_P$, there is a $\PP^2$--family of such $C_b$, with free parameter $b_5$ and $B$ varying along $l^i_P$. Similarly for $\Gamma_{\iota (U)}\bigcap \Lambda$.
  
 We now  write the equations of the conics $\Gamma_U \bigcap \Lambda$ and $\Gamma_{\iota (U)}\bigcap \Lambda$ in coordinates $[x:y:z]$ on $\Lambda.$ Let $V=[xP+yQ+zR]=[v_i]_i$, with $v_i=x+yu_i+zl_i$. The quadratic polynomial defining $\Gamma_U \bigcap \Lambda$ is
\bean  Q_U(V)&= &\sum_{i=0}^4M_i(U)(x+yu_i+zl_i)^2 \quad  \mbox{ hence }\\ \label{1st quadratic} Q_U(V) &= & y^2 (\sum_{i=0}^4M_i(U)u^2_i)+ 2yz (\sum_{i=0}^4M_i(U)u_il_i)+ z^2 (\sum_{i=0}^4M_i(U)l^2_i),  \eean
  since the missing terms came with coefficients \bea \sum_{i=0}^4M_i(U)=\sum_{i=0}^4M_i(U)u_i=\sum_{i=0}^4M_i(U)l_i =0.\eea
Similarly for  $\Gamma_{\iota (U)}\bigcap \Lambda$ we have 
 \bean  Q_{\iota (U)}(V)&= &\sum_{i=0}^4M_i(U)u_i^3(x+yu_i+zl_i)^2 \\  \label{2nd quadratic} Q_{\iota (U)}(V) &= & x^2 (\sum_{i=0}^4M_i(U)u^3_i)+ 2xz (\sum_{i=0}^4M_i(U)u^3_il_i)+ z^2 (\sum_{i=0}^4M_i(U)u^3_il^2_i),\eean   since the missing terms have coefficients  \bea \sum_{i=0}^4M_i(U)u_i^5=\sum_{i=0}^4M_i(U)u_i^4=\sum_{i=0}^4M_i(U)u^4_il_i =0 \eea
  Thus each of $\Gamma_U \bigcap \Lambda$ and $\Gamma_{\iota (U)}\bigcap \Lambda$ is a union of two lines in $\Lambda$, passing through $P$ and $Q$ respectively. 
  \epr
  
  \bl \label{discriminant}
  Consider quadratic polynomials associated to the scaled discriminants of (\ref{1st quadratic}) and (\ref{2nd quadratic}), respectively:
  \bea  \Delta(V) & = & (\sum_{i=0}^4M_i(U)u_iv_i)^2-(\sum_{i=0}^4M_i(U)u^2_i)(\sum_{i=0}^4M_i(U)v^2_i), \mbox{ and } \\ \Delta'(V)& = & (\sum_{i=0}^4M_i(U)u^3_iv_i)^2-(\sum_{i=0}^4M_i(U)u^3_i)(\sum_{i=0}^4M_i(U)u^3_iv^2_i)\eea 
  $\Delta$ and $\Delta'$ are exchanged under the existing involution (up to a factor of $e_5^{10}$), and 
  \bea \Delta(xP+yQ+zR) = \Delta(R)z^2, \quad \quad \Delta'(xP+yQ+zR) = \Delta'(R)z^2, \quad \quad \mbox{ for all } R\in \Lambda \eea
  \el 
  
  \bpr Exchange under involution (up to scalar) can be checked directly.
  Evaluating the expressions above at $V=P=[1:1:1:1:1]$, and  $V=Q=[u_i]_i$ directly gives  $\Delta(Q)=0$ and $\Delta'(P)=0$. On the other hand, as $p_i^5=M_i(U)$ form the solution in (\ref{eq Mi}),
  \bea  \sum_{i=0}^4 M_i(U) =  \sum_{i=0}^4 M_i(U)u_i = \sum_{i=0}^4 M_i(U)u_i^4 = \sum_{i=0}^4 M_i(U)u_i^5 =0\eea 
%   \bea  \Delta(V) =  -\sum_{i>j}^4M_i(U)M_j(U)(u_iv_j-u_jv_i)^2 \mbox{ and } \Delta'(V) =  -\sum_{i>j}^4M_i(U)M_j(U)u_i^3u_j^3(v_i-v_j)^2,\eea 
  hence $\Delta(P)=0$ and $\Delta'(Q)=0$. For general $v_i=x+yu_i+zl_i$, we get 
  \bea  \Delta(V) &= & -\sum_{i>j}^4M_i(U)M_j(U)\left(x(u_i-u_j) + z(u_il_j-u_jl_i)\right)^2 \\ & = & -2xz\sum_{i>j}^4M_i(U)M_j(U)(u_i-u_j)(u_il_j-u_jl_i) +
  %z^2 \sum_{i>j}^4M_i(U)M_j(U)\left(u_il_j-u_jl_i\right)^2= 
  z^2\Delta(R), \quad \mbox{ and  }
   \eea
 \bea & &  \sum_{i>j}^4M_i(U)M_j(U)(u_i-u_j)(u_il_j-u_jl_i) = \sum_{i,j}^4M_i(U)M_j(U)u_iu_jl_i-  \sum_{i,j}^4M_i(U)M_j(U)u_i^2l_j\\ & & \quad \quad \quad \quad \quad   = (\sum_{i}^4M_i(U)u_il_i)(\sum_{j}^4M_j(U)u_j)-  (\sum_{i,j}^4M_i(U)u_i^2)(\sum_{j}M_j(U)l_j)=0\eea  
  which shows that $\Delta(V)=\Delta(R)z^2$. Similarly for $\Delta'.$
  \epr
  
  The above Lemma indicates that any choice of point $R$ in $\Lambda \setminus PQ$ is equally suited for factoring the quadratic polynomials (\ref{1st quadratic}) and (\ref{2nd quadratic}).

  \begin{notation}
  Let $S_{mn}(U)$ be the symmetric polynomials given by
  \bean \label{sums} \sum_{i=0}^4M_i(U)u_i^ml_i^n = \delta(U) S_{mn}(U). \eean 
  \end{notation}
 
  These are well defined as the sums $T_{mn}(U):=\sum_{i=0}^4M_i(U)u_i^ml_i^n$ are symmetric polynomials in $\CC[u_i]_i$ which are multiples of the discriminant $\delta (U)$. Indeed, $T_{mn}(U)|_{u_j=u_k}=0$ because in the case $u_j=u_k$ we have $M_i(U)=(-1)^id_in_i=0$ for all $i\not= j, k$ while $M_j(U)=-M_k(U)$ and hence  $l_j=l_k$  by equations (\ref{eqs for lambda}) applied to $R=[l_i]_i\in \Lambda$ (in the case $M_j(U)=-M_k(U)\not=0$).
  
 Recall that $\Delta$ and $\Delta'$ are swapped under the involution (up to a factor of $e_5^{10}$). Up to a factor of $e_5^5$, the involution $\tau$ described in (\ref{involution}) swaps $S_{m,n}(U)$ with $S_{(5-m-n),n}(\iota(U))$.

   %translates to  \bean \label{involution on S}  S_{m,n}(U)= S_{(5-m-n),n}(\iota(U)). \eean
 %  \br times something from discriminant and then careful about the rule for e-s \er
   
     Below are the relevant values of $S_{m,0}$ in terms of elementary symmetric polynomials: 
  \bean  \label{first Ss} \begin{array}{llllll}  S_{00}=0, & S_{10}=0, & S_{20}=e_2, & S_{30}=e_3, & S_{40}=0, & S_{50}=0. \end{array}\eean 
  
  Indeed, the null terms are due to equation (\ref{eq Mi}) with solutions $p_i^5=M_i(U)$, where $(-1)^iM_i(U)$ are the $4\times 4$ minors of the associated matrix. By the same token, 
   \bea T_{20} =   \left|\begin{array}{ccccc} 1&1&1&1&1 \\
                          u_0 & u_1 & u_2 & u_3 & u_4 \\
                           u_0^2 & u_1^2 & u_2^2 & u_3^2 & u_4^2 \\
                          u_0^4 & u_1^4 & u_2^4 & u_3^4 & u_4^4 \\
                          u_0^5 & u_1^5 & u_2^5 & u_3^5 & u_4^5 
 \end{array} \right|  \quad \mbox{ and } \quad T_{30}  =  \left|\begin{array}{ccccc} 1&1&1&1&1 \\
                          u_0 & u_1 & u_2 & u_3 & u_4 \\
                           u_0^3 & u_1^3 & u_2^3 & u_3^3 & u_4^3 \\
                          u_0^4 & u_1^4 & u_2^4 & u_3^4 & u_4^4 \\
                          u_0^5 & u_1^5 & u_2^5 & u_3^5 & u_4^5 
 \end{array} \right|  \eea 
are  $5\times 5$ minors of a $6\times 6$ Vandermonde matrix $V=(u_i^{j})_{0\leq i, j \leq 5}$. We note that $T_{23}$ and $-T_{30}$ are entries in the adjoint $V_{\mbox{adj}}$ of the Vandermonde matrix, which can be calculated based on the matrix multiplication formula $V_{\mbox{adj}}V=\prod_{0\leq i<j\leq 5}(u_j-u_i)\textbf{I}_6.$ Hence $T_{23}$ and $-T_{30}$ are
the coefficients of $u^3$ and $u^2$ in the polynomial 
\bea  p(u)= \delta(U)\prod_{i=0}^4(u-u_i)= \delta(U)(u^5-e_1 u^4+ e_2 u^3 -e_3 u^2 +e_4 u -e_5). \eea 
Consider the extension of algebras $\CC[u_i]_i \hookrightarrow \CC[u_i]_i[\sqrt{\Delta(R)}, \sqrt{\Delta'(R)}].$ From now on, we will work in the associated $4$-to-$1$ cover of $\PP^4.$
  
\begin{corollary} \label{4 pencils}
 Consider the $\PP^5$--family of conics passing through $P$ and $Q$ given by parametric equations (\ref{conic parametrisation}). Let $B:=[b_0:b_1:b_2:b_3:b_4]$. Then within this family there are in general four $\PP^1$-- families of conics  which satisfy
\bea  C_b \subset \Lambda \mbox{  and } C_b\cdot X_{\Lambda}\geq 3P + 3Q \mbox{ in } \Lambda \eea
with $b=[B:b_5 ]$ where $B=[\alpha_{\pm}+\beta_{\pm}u_i+l_i]_{i\in\{0,1,2,3,4\} }$ and $[\alpha_{\pm}:1]$ are the solutions of the quadratic equation (\ref{1st quadratic}) and $[\beta_{\pm}:1]$ are the solutions of the quadratic equation (\ref{2nd quadratic}) for a suitable choice of point $R=[l_i]_i\in \Lambda.$ Namely, if $S_{20}\not=0$ and $S_{30}\not=0$ then
\bea S_{20} \beta_{\pm}= - S_{11}\pm \sqrt{\Delta(R)}, \quad \quad S_{30} \alpha_{\pm}= - S_{31}\pm \sqrt{\Delta'(R)},  \eea 
$\Delta(R)=S_{11}^2-S_{20}S_{02}$ and $\Delta'(R)=S_{31}^2-S_{30}S_{32}.$

% the case when $S_{20}=0$ then $\beta = 0$
\end{corollary}

Thus we get 4 families of conics exactly when $S_{20}\not=0$, $S_{30}\not=0$  and $\Delta(R)\not=0$, $\Delta'(R)\not=0.$
 
 \section{Equations for the conic and cubic curves in suitable singular plane sections of the Fermat quintic.} \label{section on conic and cubic}
 
Assume that for a chosen plane $\Lambda$ with marked points $P$ and $Q$, the section  $X_\Lambda$ decomposes into a conic $C$ and a cubic $C'$. In this section we will calculate the equations of $C$ and $C'$ in $\Lambda$, in terms of the coordinates $[x:y:z]$ on $\Lambda$ determined by the coordinate simplex $PQR.$
  
 The conic $C$ will be a member in one of the pencils described in Corollary \ref{4 pencils}. Let this pencil be parametrised by $b=[B:b_5 ]$ where $B=\alpha P+ \beta Q + R = (\alpha+\beta u_i+l_i)_{i\in\{0,1,2,3,4\} }$ for a choice of roots $\alpha$ and $\beta$ for the quadratics equations in (\ref{1st quadratic}) and (\ref{2nd quadratic}) respectively. For simplicity we will denote $\lambda=b_5^2$.
 %thus working with $[B:c]\in\PP[1:1:1:1:2]$. 
 
 In the chosen coordinates on $\Lambda$ for which $P=[1:0:0]$, $Q=[0:1:0]$ and $R=[0:0:1]$, the parametric equations (\ref{conic parametrisation}) translate to 
  \bea \label{plane conic parametrisation} \gamma_{b}: \PP^1 \ra \PP^2, \quad \quad \quad \quad \gamma_{b}([t:s])= [t^2 P+ ts B + \lambda s^2 Q]=[t^2+\alpha ts: \beta ts +\lambda s^2: ts ]. \eea 
  After the change of coordinates on $\Lambda\cong \PP^2$
  \bea  \quad \quad \quad \quad  [x: y: z] \ra \left[ x-\alpha z: z: \frac{1}{\lambda }(y-\beta z) \right]\eea
  the parametrisation of $C$ becomes the Veronese embedding, hence $C$ has defining polynomial
  \bean  \label{conic equation}  f(x,y,z)= (x-\alpha z)(y-\beta z) -  \lambda z^2  = (\alpha\beta - \lambda ) z^2 -(\beta x+ \alpha y) z + xy.\eean 
  Thinking of this as a polynomial in $z,$ we denote the coefficients 
  \bea f_2=d:= \alpha\beta - \lambda, \quad \quad f_1=-(\beta x+ \alpha y), \quad \quad f_0= xy.\eea 

\subsection{Equations for the paired cubic}  Assume $fg=F_{\Lambda}$ where $F_{\Lambda}$ is the polynomial of the quintic section $X_{\Lambda} $. We will write the cubic polynomial $g(x,y,z)=\sum_{k=0}^3 g_kz^k$, with $g_k(x,y)$ homogeneous polynomial of degree $(3-k)$. 
  We also write the equation of the quintic section $X_\Lambda$, noting the common factor $\delta(U)$ which comes from (\ref{sums}) above. 
  
  \bea  \delta(U) F_\Lambda(x,y,z)=\sum_{i=0}^4 M_i(U) (x+ u_i y+ l_i z)^5 = \delta(U) \sum_{j+l+k=5}\frac{5!}{j ! l ! k !} S_{jk}(U)x^ly^jz^k\eea 
  and based on the initial values of $S_{jk}(U)$ as listed in (\ref{first Ss}) above.

   \bea F_\Lambda(x,y,z) & = &  10x^2y^2(S_{20}x+S_{30}y)+(20x^3yS_{11} + 30 x^2y^2 S_{21}+20xy^3S_{31})z +\\ & + & \sum_{k=2}^4 \sum_{j+l=5-k}\frac{5!}{j ! l ! k !} S_{jk}x^ly^jz^k  \eea 
  
  \begin{notation}
    We denote $k(x,y):=S_{20}x+S_{30}y = e_2x+e_3y$ and $l(x,y):= \beta x + \alpha y.$  \end{notation}
    
By identifying coefficients $\sum_{i+j=k}g_if_j=F_k$ for $k\in\{0,1,2,3\}$,
%  \bea g_0f_0=F_0, \quad  g_1f_0+ g_0f_1= F_1, \quad g_2f_0+ g_1f_1 +g_0f_2= F_2, \quad  g_3f_0+ g_2f_1+g_1f_2+g_2f_1+g_3f_0= F_3 \eea 
we determine the coefficients $g_k(x,y)$ of the cubic equation recursively: 
  \bean \label{g0}  g_0(x,y)& = & 10xy(S_{20}x+S_{30}y)=10xy k(x,y); \\ 
 \label{g1}  g_1(x,y) & = &(20x^2S_{11} + 30 xy S_{21} + 20y^2S_{31}) + 10k(x,y)l(x,y);\\ 
 \label{g2}  xyg_2 & = & F_2+ g_1 l - g_0 d; \\
 \label{g3}  xyg_3 & = & F_3+ g_2 l - g_1 d.\eean
  Differentiating equation (\ref{g2}) with respect to $x$ and $y$, and since $g_2=xg_{2 x}+yg_{2 y}$ we get
  %as $\mbox{deg}(g_2)=1$, we obtain
  \bean \label{g2 final} g_2  = \frac{1}{2}(F_{2 xy}+ \alpha g_{1 x}+ \beta g_{1 y}+l g_{1 xy}- d g_{0 xy}); \eean
  Here $F_2(x,y)=10(S_{02}x^3+3S_{12}x^2y+3S_{22}xy^2+S_{32}y^3)$ hence $F_{2 xy}(x,y)=60(S_{12}x+S_{22}y).$ Also $g_{0 xy}(x,y)=20k(x,y)$ and $g_{1 xy}=30S_{21}+10k(\alpha,\beta)$, while $g_{1 y}, g_{1 x},$ are of degree 1 hence evaluating at $P=(1,0)$ and $Q=(0,1)$ we may rewrite
  \bea g_{1 x}&=& x g_{1 xx}+y g_{1 xy}=2xg_1(P)+yg_{1 xy}, \\g_{1 y}&=&x g_{1 xy}+y g_{1 yy}=xg_{1 xy}+2yg_1(Q).\eea 
%  This leads to the final formula for $g_2$:
  Similarly for the constant coefficient $g_3$ we get
  \bea  g_3  = F_{3 xy}+ \alpha g_{2 x}+ \beta g_{2 y}- d g_{1 xy}. \eea
Here  $F_3(x,y)=10(S_{03}x^2+2S_{13}xy+S_{23}y^2)$ hence $F_{3 xy}(x,y)=20S_{13}.$ Also since  $\mbox{deg}(g_2)=1,$ 
%then evaluating at $P$ and $Q$ we may rewrite
% \bea \alpha g_{2 x}+ \beta g_{2 y}=g_2(\alpha, \beta)=30(\alpha S_{12}+\beta S_{22})+ g_1(\alpha, \beta) - 10 d k (\alpha, \beta),\eea
%(having applied equation (\ref{g2 final}) for $g_2$ and Euler's relation for $g_1$). Also calculating $g_{1 xy}$, 
 \bea g_{2 x}= g_{2}(P) \quad  \mbox{ and } g_{2 y}= g_{2}(Q).\eea

%\bean \label{g3 final} g_3  = 20S_{13}+ 30 S_{12} \alpha +30 S_{22} \beta + g_1(\alpha, \beta) - d(20k (\alpha, \beta)+30S_{21}). \eean 

 The $[x:y]$ coordinates of the intersection points of $C$ and $C'$ are given by the zeroes of the resultant $R(f,g,z)$ which is a degree 6 homogeneous polynomial in $x,y$, written in terms  of $(f_i)_i$ and $(g_i)_i$. We know that $R(f,g,z)$ is a multiple of $xy$ since $P, Q \in C\bigcap C'$.

\section{Equations for the the parameter space of conics in the Fermat quintic.} \label{section on eqs}

We continue with the set-up where we parametrise plane sections together with two marked singular points. Previously we wrote the equations of the conic and cubic components of a plane section $X_{\Lambda},$ which were obtained from identifying coefficients of $z^k$ in the equation of $X_{\Lambda},$ for $k\in\{0,1,2,3\}$.
The equations for $k\in\{3, 4,5\}$ will yield 5 necessary and sufficient conditions for the existence of such conics in $X$:
The equation (\ref{g3}) requires that $xy|(F_3+ g_2 l - g_1 d)$ ( the equivalent conditions from  equations (\ref{g0})--(\ref{g2}) are automatically satisfied due to the context set-up in the previous section). Equivalently, the polynomial $F_3+ g_2 l - g_1 d$ is zero when evaluated at $P=(1,0)$ and $Q=(0,1)$. Similarly for  $k\in\{ 4,5\}$. We obtain

\bt \label{main equations} Let $\Lambda$ be a plane in $\PP^4$ spanned by the points $P=[1:1:1:1:1]$, $Q=[u_0:u_1:u_2:u_3:u_4]$ and $R=[l_i]_i$. Consider the sums $S_{mn}$ defined by equation (\ref{sums}), and the extension $\CC[u_0,u_1,u_2,u_3,u_4][\alpha,\beta]$ with $\alpha$, $\beta$ satisfying 
\bea  S_{30}\alpha^2 +2S_{31}\alpha+S_{32}   \quad \mbox{ and } \quad S_{20}\beta^2 +2S_{11}\beta+S_{02}. \eea Then the plane spanned by $P, Q, R$ cuts the Fermat quintic along a conic and a cubic (possibly degenerate) passing through $P$ and $Q$ iff the following equations are satisfied for some parameter $d$:
\bean \label{eq 1} 10 S_{03} + \beta g_2(P)  -  g_1(P) d & = & 0 \quad  \quad  \mbox{ and } \\  \label{eq 2} 10 S_{23} + \alpha g_2(Q) -  g_1(Q) d  & = & 0\\ \label{eq 3} 5 S_{04} + \beta g_3  -  g_2(P) d & = & 0 \quad  \quad  
\\ \label{eq 4} 5 S_{14} + \alpha g_3 -  g_2(Q) d  & = & 0\quad  \quad  \mbox{ and } \\ \label{eq 5}  S_{05} - g_3 d & = & 0.
\eean 

\et

Equations (\ref{eq 1}) and (\ref{eq 2}) above are linear in $d$; the remaining equations are quadratic  in $d$, since both $g_2$ and $g_3$ are linear in $d$. In detail, the coefficients are
\bean \label{first eq details} && g_1(P) = 20 S_{11}+ 10 S_{20} \beta, \quad \quad g_1(Q)=20 S_{31}+ 10 S_{30} \alpha, \\
&& g_2(P) = 30 S_{12}+\alpha g_1(P)+ \beta g_{1 xy}- 10 S_{20}d \\
&& g_2(Q) = 30 S_{22} + \beta g_1(Q) +\alpha g_{1 xy}-10 S_{30}d, \\
%&& g_2(P) = 30 S_{12}+20  S_{11} \alpha + 15 S_{21} \beta +15 S_{20}\alpha \beta +5S_{30}\beta^2 - 10 S_{20}d, \\ \label{last eq details}
%&& g_2(Q) = 30 S_{22}+20  S_{31} \beta + 15 S_{21} \alpha+15 S_{30}\alpha \beta +5S_{20}\alpha^2 - 10 S_{30}d
&&  \label{last eq details} g_3  = 20S_{13}+ \alpha g_{2}(P)+ \beta g_{2}(Q)- d g_{1 xy}, \eean 
where as seen previously $g_{1 xy}=30S_{21}+10k(\alpha,\beta)=30S_{21}+10S_{20}\alpha+10S_{30}\beta.$
%and $g_3$ as given in formula (\ref{g3 final}).
%\br triple-check these; discuss the symmetry through tau. can we use it somehow to argue that the equations are dependent? also check out relations among Smn \er   

Equation $(\ref{eq 5})$ is invariant under the involution $\tau$, while the same involution swaps (\ref{eq 1}) and (\ref{eq 2}) as well as  (\ref{eq 3}) and (\ref{eq 4}). 
Hence for symmetry reasons it seems optimal to substitute $d$ from (\ref{eq 1}) into (\ref{eq 3})  and $d$ from (\ref{eq 2}) into (\ref{eq 4}) thus obtaining two equations in $\CC[u_i]_i[\alpha, \beta]$ which are pairs under $\tau$. We can use a combination of $(\ref{eq 1})$ and $(\ref{eq 2})$ to eliminate $d$ from (\ref{eq 5}). Finally, $(\ref{eq 1})$ and $(\ref{eq 2})$ together produce a fourth equation in $\CC[u_i]_i[\alpha, \beta]$.

%while also $\mbox{deg}(g_1)=2$ hence $2g_1=xg_{1 x}+yg_{1 y}$ and $2g_1(P)=g_{1 x}(P)$, $2g_1(Q)=g_{1 y}(Q)$ so 

\subsection{Choices of coordinate simplex for the plane $\Lambda$} 

The equations for the conics, cubics in $X$ and their parameter spaces found in the previous section all depend on the choice of coordinates the plane $\Lambda.$ The coordinate simplex we worked with has two fixed vertices $P$ and $Q$, while the third vertex $R$ can be chosen to be any point on $\Lambda \setminus PQ.$
This choice will determine the format of the entries $S_{mn}$ in the equations above, the format of $\alpha, \beta$ and the parameter $d.$ At the same time, the degree 4 extension  $\CC[u_i]_i[\alpha, \beta]$  of $\CC[u_i]_i$ does not depend on the choice of the point $R$ due to Corollary \ref{4 pencils}.

 It will be interesting to determine exactly how the system of equations found in the previous section varies with the choice of $R$, and use this to gain insights into the geometry of the space of conics determined by these equations. Within the constraints of this article we were not able to delve deeper into this topic but we hope to do this elsewhere.  

In this section we will only briefly discuss and compare two particular choices for the point $R\in \Lambda$, their advantages and drawbacks. Recall that the plane $\Lambda$ was given by the equations (\ref{eqs for lambda}). Recall that $\delta(U)S_{mn}=\sum_{i=0}^4M_i(U)u_i^ml_i^n.$

\textbf{Option 1.} Let $R=[l_i]_i$ be given by the equations
\bea l_i = S_{8 0}u_i^5-S_{9 0}u_i^4 \eea 
Indeed, $R$ satisfies the equations (\ref{eqs for lambda'}) as we have \bea \sum_{i=0}^4M_i(U)u_i^4l_i&=& S_{8 0}\sum_{i=0}^4M_i(U)u_i^9-S_{9 0}\sum_{i=0}^4M_i(U)u_i^8 =S_{8 0}S_{9 0}-S_{9 0}S_{8 0}=0, \\ \sum_{i=0}^4M_i(U)l_i&=& S_{8 0}\sum_{i=0}^4M_i(U)u_i^5-S_{9 0}\sum_{i=0}^4M_i(U)u_i^4 =S_{8 0}S_{5 0}-S_{9 0}S_{40}=0 \eea as $S_{50}=S_{40}=0$ from the construction of $M_i(U)$ in the proof of Proposition \ref{double tangents}.

The terms $S_{m0}$ have initial values: 
  \bea  \begin{array}{llllll}  S_{00}=0, & S_{10}=0, & S_{20}=e_2, & S_{30}=e_3, & S_{40}=0, & S_{50}=0. \end{array}\eea 
  All further values can be found recursively from the Newton-type identity \bea  S_{m,0}=e_1S_{m-1,0}-e_2S_{m-2,0}+ e_3S_{m-3,0}-e_4S_{m-4,0}+e_5S_{m-5,0},\eea
  which follows from the identity $\sum_{i=0}^4M_i(U)q(u_i)u_i^{m-5}=0$ for \bea q(u):=\prod_{i=0}^4(u-u_i)=u^5-e_1u^4+e_2u^3-e_3u^2+e_4u-e_5.\eea
In turn, the general terms $S_{mn}$ are determined by the recurrence
\bea  S_{mn}= \sum_{k=0}^n (-1)^{n-k} \binom{n}{k} S_{80}^k S_{90}^{n-k}S_{m+4n+k,0}. \eea 
Indeed, this follows from applying the binomial formula to the term $l_i^n =(S_{8 0}u_i^5-S_{9 0}u_i^4)^n$ in the definition of $S_{mn}.$

We thus have an algorithmic procedure of finding all the coefficients entering the equations of the space of conics. However, the complexity of the terms  $S_{mn}$ as polynomials of $e_i$ increases fast with $m,n.$ For the time being we have been unable to extract geometrically meaningful information from this procedure, but it's possible that more patterns will emerge from the interplay of the algebra and geometry here. 
%For example, the invariance of the quadratic forms $\Delta, \Delta'$ will allow us to compare. 

\textbf{Option 2.} 
Our second choice $R'$ is defined geometrically as the point of intersection of the planes $\Lambda=\Lambda(U)$ and $\Lambda'$ 
%$\Lambda'=\Lambda(\iota(U))$ 
in $\PP^4$, where $\Lambda'$ is given by equations
    \bean \label{eqs for lambda'}  \sum_{i=0}^4M_i(U)u_i^2v_i = 0 \quad \quad \mbox{ and } \quad \quad \sum_{i=0}^4M_i(U)u^3_iv_i = 0, \eean
     
  With the notations from the previous sections, up to a $\CC^*$-- factor we have $R'=[l'_i]_i$ with 
  \bean \label{third point} l'_i= e_2(i) \prod_{j\not= i} n_j, \eean
  where $n_j$ are as defined in Theorem \ref{P4 parametrisation}. %Thus $l'_i$ are symmetric polynomials in $(u_j)_{j\not= i}.$
  
  Below are the first few values of $S'_{m,n}$ in terms of elementary symmetric polynomials: 
  \bea  \begin{array}{llllll}  S'_{00}=0, & S'_{10}=0, & S'_{20}=e_2, & S'_{30}=e_3, & S'_{40}=0, & S'_{50}=0, \\  S'_{01}=0, & S'_{11}=0, & S'_{21}=\prod_{i}n_i, & S'_{31}=0, & S'_{41}=0.  \end{array}\eea 
  Comparing with the coefficients in formulae (\ref{first eq details} - \ref{last eq details}), we see that simplifies the system (\ref{eq 1}--\ref{eq 5}) and in particular eliminates the parameter $d$. 
All the remaining terms are products of $\prod_{i}n_i$ due to our choice in formula (\ref{third point}). Then by computer assisted calculations we verified \bea S'_{02}&=& -e_2e_3^2S_{80}\prod_{i}n_i, \quad \quad \\
%  \bea  S'_{02} & = & -e_2e_3^2(e_3^2e_1^2-e_3^2e_2-e_4e_2e_1^2+e_4e_2^2-e_4e_3e_1+e_5e_2e_1+e_5e_3)\prod_{i}n_i= -e_2e_3^2S_{80}\prod_{i}n_i.\\  S'_{12} & = & -e_3^2(e_3^2e_2^2e_1-e_3^3e_1^2-e_4e_2^3e_1+e_4e_3e_2e_1^2-e_4e_3e_2^2+e_4e_3^2e_1+e_5e_2^3-e_5e_3e_2e_1-e_5e_3^2)\prod_{i}n_i. \eea
S'_{12}e_2+S'_{02}e_3&=&-e_2e_3^2(e_3^2e_2^2e_1+e_5e_3^3+e_4^2e_3e_2e_1-e_3^2e_2)\prod_{i}n_i. \eea
The next term $S'_{22}$ is the transform of $S'_{12}$ under the involution $\tau$ and $S'_{32}$ is the transform of $S'_{02}$ under the same involution (up to a factor of $e^k_5$).
%related by formula (\ref{involution}), namely  $S'_{22}(U)=S'_{12}(\iota(U))e_5^4(\iota(U))$, where $S'_{12}(\iota(U))$ is calculated based on $e_j(\iota(U))=e_{5-j}(U)/e_5(U)$. Similarly $S'_{32}$ is the pair of $S'_{02}$:
%\bea  S'_{22} & = & -e_2^2(e_4e_3^2e_2^2-e_4e_3^3e_1-e_4^2e_2^3+e_4^2e_3e_2e_1-e_5e_3^2e_2e_1+e_5e_3^3+e_5e_4e_2^2e_1-e_5e_4e_3e_2-e_5^2e_2^2)\prod_{i}n_i; \\
%  S'_{32} & = & -e_2^2e_3(e_4^2e_2^2-e_4^2e_3e_1-e_5e_3e_2^2+e_5e_3^2e_1-e_5e_4e_2e_1+e_5e_4e_3+e_5^2e_2)\prod_{i}n_i.   \eea
%We note that the discriminants in this case are $\Delta(R')=-S'_{02}S'_{20}$ and $\Delta'(R')=-S'_{32}S'_{30}$.

The full pattern of identities among the terms $S_{mn}$ and their significance for the solutions of the equations in the previous section is yet to emerge. 

\section{Plane sections in the exceptional locus of $\overline{U}_2(X) \to \overline{\Delta}_{2,2}(X).$}

In this section we will illustrate the application of the equations obtained in Theorem \ref{main equations} to the study of conics found for plane sections in some of the components of the exceptional locus of $\overline{U}_2(X) \to \overline{\Delta}_{2,2}(X).$ 

While the equations in Theorem \ref{main equations} were  set up on 
%the 4-to-1 cover of $\PP^4$ which results from combining Theorem \ref{double tangents} with Proposition \ref{4 pencils}, 
$\Delta_{2,2}(X)\setminus X_{\sigma}$, we note that these equations can be extended to work for points outside the regular sets of maps $\phi$ and $m$. Namely, for a triple $(\Lambda, P, Q)\in \overline{U}_2(X)$, we can unpack the change of coordinates $\nu_P^{-1}$ introduced in Section 4, so that the equations (\ref{first eq details} -- \ref{last eq details}) can be translated following the rules
\bean \label{new Smn} M_i(U)=p_i^5, \quad  q_i=p_iu_i, \quad  h_i=p_il_i \quad \mbox{ and hence } \quad  S_{mn} =\sum_{i=0}^4p_i^{5-n-m}q_i^mh_i^n. \eean 
Here $h_i$ are the coordinates of the point $R$ introduced in Section 4.
Rather than reversing the benefits of Theorem \ref{double tangents}, this coordinate choice can be interpreted as working over a blow-up of $\PP^4$ as discussed in section \ref{resolution} (an intermediate space between  $\PP^4$ and $\overline{M}_{0,7}$ like the Losev-Manin spaces would be sufficient). 

Consider the case when $(\Lambda, P, Q)$ is in the exceptional locus of $f: \overline{U}_2(X) \to \overline{\Delta}_{2,2}(X)$. Recall that this is the locus in the preimage of  $X_{[k]}=\sigma_{[k]}(\Delta_{\PP^4})\bigcap (X\times X)$ for the action defined by formula (\ref{Z4 locus }), namely
$P=[p_i]_i$ and $Q=[p_i\sqrt{-1}^{k_i}]_i$ with $[k_i]_i\in \ZZ_4^4$. 

Recall that in this case, $\overline{\cT}_{X,P}=\overline{\cT}_{X,Q}$ and hence $\Lambda$ varies in the pencil given by $PQ \subset \Lambda \subset \overline{\cT}_{X,P}.$ As well, $P$ and $Q$ are both singular points of the quintic surface $X\bigcap \overline{\cT}_{X,P}=V(\sum_{i=0}^4 x_i^5, \sum_{i=0}^4 x_i p_i^4).$ The tangent cones $C_P(X)$ and $C_Q(X)$ to the surface $X\bigcap \overline{\cT}_{X,P}$ at $P$ and $Q$ respectively, are quadratic surfaces  defined in $\overline{\cT}_{X,P}$ by the polynomials
\bea  \sum_{i=0}^4 x_i^2 p_i^3 \quad  \mbox{ and }  \quad \sum_{i=0}^4 x_i^2 q_i^3=\sum_{i=0}^4 x_i^2 p_i^3\sqrt{-1}^{3k_i} \eea
respectively. These correspond to the quadratic equations in Proposition \ref{conic parametrisation}.

To parametrise $\Lambda$ as above, we will pick two distinct points $S, T\in \overline{\cT}_{X,P} \setminus PQ$ so that together $P, Q, S, T$ span $\overline{\cT}_{X,P},$ and we will let $R=uS+vT=[h_i]_i$ for $[u:v]\in \PP^1.$ We can then define $\Lambda$ as the projective plane spanned by $P,Q, R.$ The tangent cones satisfy $X_\Lambda \bigcap C_P(X)= C_P(X_\Lambda)$ and $X_\Lambda \bigcap C_Q(X)= C_Q(X_\Lambda)$ and Corollary \ref{4 pencils} applies, with 
%quadratic equations (\ref{1st quadratic}) and (\ref{2nd quadratic}) taking the form \bean \label{1st new quadr} Q_{C_P(X_\Lambda)}&= & y^2 (\sum_{i=0}^4p_i^3q^2_i)+ 2yz (\sum_{i=0}^4p_i^3q_ik_i)+ z^2 (\sum_{i=0}^4p_i^3k^2_i)  \\ \label{2nd new quadr} Q_{C_Q(X_\Lambda)}&= & x^2 (\sum_{i=0}^4p_i^2q^3_i)+ 2xz (\sum_{i=0}^4p_iq_i^3k_i)+ z^2 (\sum_{i=0}^4q^3_ik^2_i).\eean  
with the notations set in (\ref{new Smn}). Then the program from sections \ref{section on conic and cubic} and \ref{section on eqs} follows accordingly.

We will illustrate how this works in the case of  $k=[0:0:0:2:2]\in \ZZ_4^4.$ We denote by $\Sigma_{[k]}$ the image of the first projection of $X_{[k]}$ on $X$, namely $\Sigma_{[k]}=V(\sum_{i=0}^4 x_i^5, \sum_{i=0}^4 x_i^5\sqrt{-1}^{k_i})$, a degree 25 surface in $\PP^4.$ For $k$ as chosen above and $\mu$ a fifth root of 1, we have \bea  \Sigma_{[k]}=\bigcup_{h\in \ZZ_5}V(\sum_{i=0}^2 x_i^5, x_3+\mu^hx_4), \eea
a union of 5 of the 50 cones making up the variety covered by lines in $X.$ We will focus on the case when $h=0$. Of course, the other cones are obtained from $X_k$ by the action of $S_5\times \ZZ_5^4,$ hence they will behave similarly in respect to conics. 

We will set the following (affine) coordinates, with $a^5+b^5+c^5=0$:
\bea && P = [a:b:c:1:-1], \quad \quad  Q = [a:b:c:-1:1], \\ &&  S = [0:-b:0:0: b^5], \quad \quad T = [0:0:-c:c^5:0].\eea 
It can be readily verified that these four points span 
$\overline{\cT}_{X,P}=V(a^4x_0+b^4x_1+c^4x_2+x_3+x_4)$ in the case $bc\not=0.$ The other cases can be treated similarly after a permutation of coordinates. For simplicity we will use an affine parameter, writing $R=tT+S$ with $t\in \AA^1$. Then since $p_i^2=q_i^2,$ the coefficients $S_{m,n}$ satisfy $S_{m0}=p_0^5+p_1^5+p_2^5+(-1)^m(p_3^5+p_4^5)=0$ for all $m\in\{0,1,2,3,4,5\}$, and 
\bea S_{01} &= & S_{21} = S_{41} = 0.\\
S_{11}  &= &  S_{31} = -2(tc^5+b^5); \\
S_{02}  &= & S_{22}=t^2(c^5+c^{10})+b^5-b^{10}; \\
S_{12} &= & S_{32}= t^2(c^5-c^{10})+b^5+b^{10};\\
S_{03}  &= & S_{23}= -t^3(c^5-c^{15})-b^5+b^{15};\\
S_{13} &= & -t^3(c^5+c^{15})-b^5-b^{15};\\
S_{14}  &= & t^4(c^5-c^{20})+b^5+b^{20};\\
S_{04} &= & t^4(c^5+c^{20})+b^5-b^{20};\\
S_{05}  &= & -t^5(c^5-c^{25})-b^5+b^{25}.
 \eea 
 In particular, the forms in (\ref{1st quadratic}) and (\ref{2nd quadratic}) equal $y$ (respectively $x$) times linear factor and admit common zeros along the line $PQ=V(z)$ (with the notations from (\ref{1st quadratic}) and (\ref{2nd quadratic})). Conics tangent to this line at both $P$ and $Q$ are those unions of lines $PQ\bigcup l \subset X$. Assuming $S_{11}\not=0$, the remaining case is when $\alpha =-S_{12}/(2S_{11})$ and $\beta =-S_{02}/(2S_{11}).$
%\br D Testa's example when $S_{11} =0?$\er 

 After simplifying and reverting to the coefficient $\lambda = \alpha\beta -d$, the equations from Theorem \ref{main equations} become 
 \bea \begin{array}{llllll} q_1 & := & 10S_{03} +30  \beta S_{12} & = & - 20 S_{11} \lambda; \\
 q_2 & := & 10S_{23} +30\alpha S_{22} & = &  - 20 S_{31} \lambda; \\
q_3& := &5S_{04}+20\beta S_{13}+20\beta ^2S_{02} & = &  - 20 S_{32} \lambda;\\
q_4& := &5S_{14}+20\alpha S_{13}+20\alpha ^2S_{12}& = &  - 20 S_{02} \lambda;\\
q_5& := & S_{05}-5S_{04}\alpha -5S_{14}\beta -60\alpha \beta S_{13}-40\alpha^2\beta S_{32}-40\alpha \beta^2S_{02} & = &  - 20 S_{13} \lambda. \end{array}
\eea
Based on the formulae for $\alpha, \beta$ and $S_{mn}$ above, the equations for $q_1$ and $q_2$ are in fact identical. Eliminating $\lambda$ and denominators from the remaining equations and substituting $b^5+c^5=-a^5$ we obtain the following polynomial equations:
\bea
f_3 & := & q_3S_{11}^2-q_1S_{12}S_{11} = -20 t^2(t-1)^2 b^5c^5a^5 =0; \\
f_4 & := & q_4S_{11}^2-q_1 S_{02}S_{11}  = -20 t^2(t-1)^2 b^5c^5a^5 =0; \\
f_5 & := & q_5S_{11}^3-q_1 S_{13} S_{11}^2 = 0,
\eea
The solution set is found based on evaluating
\bea   \begin{array}{ll} f_5|_{t=0}=2 b^{20}(b^{20}-1), & f_5|_{t=1}=2a^4(a^{20}-1),  \\  f_5|_{c=0}=2b^{20}(b^{20}-1), & f_5|_{b+c=0}=2(t-1)^8b^{20}(b^{20}-1)
\end{array} \eea and similarly for $b=0$ and $b+\mu^hc=0,$ where $\mu$ is a fifth root of unity.

%and $X_{\sigma}:=\ZZ_4^4\cdot(\Delta_{\PP^4})\bigcap (X\times X)=\bigcup_{[k]\in\ZZ_4^4}X_{[k]}$.

The solutions split into 3 symmetric cases: $t=0$ and $b^{20}(b^{20}-1)=0$, or $t=\infty$ and $c^{20}(c^{20}-1)=0$ or $t=1$ and $a^4(a^{20}-1)=0$. In the last case let $a^4=\mu.$ Then the points $P, Q$ and $R$ satisfy the equations
\bea c(x_0+\mu^{-1}x_3+\mu^{-1}x_4)= ax_2 \quad \mbox{ and } \quad cx_1=bx_2, \eea 
Thus the solutions $\Lambda$ are in the orbit of the family in Example \ref{S3 conics} under the automorphism group $S_5\times \ZZ_5^4.$

It remains to look at the case when $S_{11}=0$, namely $t=-b^5/c^5.$ Then either our conic is the double line $PQ$ or $S_{12}=S_{02}=0$ hence $c^{5}=-b^{5}$ and $t= 1.$ The solutions in this case are contained in the orbit of the family in Example \ref{Z2 symmetries}.
Perhaps not surprisingly, we obtain the following

\begin{corollary}
 All the planes containing both a line and a smooth conic in $X$ are in the orbits of Examples \ref{S3 conics} and \ref{Z2 symmetries} under the natural action of the automorphism group $S_5\times \ZZ_5^4.$ 
\end{corollary}

This is not exhaustive study of the consequences of Theorem \ref{main equations}. A full systematic study of the solution set of the resulting equations, as well as of the exceptional loci involved in our constructions will require a more extensive discussion which is yet to be finalised. In the general case, it is likely that more insights will be needed to find patterns and reduce the complexity of the equations. Optimistically, we would hope that these will be able to be brought to a point where the degree 2 equivalent of equation \ref{contribution to lines } can be verified. 

The author is grateful to the referee for very useful feedback.

\providecommand{\bysame}{\leavevmode\hbox to3em{\hrulefill}\thinspace}

\end{document}